\documentclass[12pt]{article}
\usepackage{graphicx}
\usepackage{amssymb, cite}
\usepackage{setspace} 
\usepackage[hmargin=1.15in,vmargin=1.15in]{geometry} 

\newcommand{\be}{\begin{equation}}
\newcommand{\ee}{\end{equation}}
\newcommand{\bqn}{\begin{eqnarray}}

\newcommand{\eqn}{\end{eqnarray}}

\newcommand{\bd}{\begin{description}}
\newcommand{\ed}{\end{description}}

\newtheorem{stat}{}[section]

\def\bs{\begin{stat}}
\def\es{\end{stat}}
\def\ben{\begin{enumerate}}
\def\een{\end{enumerate}}

\def\bp{\noindent{\bf Proof}  \ \ \ }

\newcommand{\ep}{\hfill $\square$}

\makeatletter
\@addtoreset{equation}{section}

\makeatother

\begin{document}

\begin{center}
{\large {\bf PACKING 3-VERTEX PATHS }}
\\[1ex]
{\large {\bf IN CUBIC 3-CONNECTED GRAPHS}}
\\[4ex]
{\large {\bf Alexander Kelmans}}
\\[2ex]
{\bf University of Puerto Rico, San Juan, Puerto Rico}
\\[0.5ex]
{\bf Rutgers University, New Brunswick, New Jersey}
\\[3ex]
{\bf The paper was accepted for publication}
\\[0.5ex]
{\bf  in Discrete Mathematics in August 14, 2009}
\end{center}

\begin{abstract}
A subgraph (a spanning subgraph) of a graph $G$ whose all components are 3-vertex paths is called an 
$\Lambda $-packing (respectively, an $\Lambda $-factor of $G$.
Let $v(G)$ and $\lambda (G)$ be the number 
of vertices and the maximum number of disjoint
3-vertex paths in $G$, respectively. 
We discuss the following old
\\[0.5ex]
\indent
{\bf Problem} (A. Kelmans, 1984).
{\em Is the following claim true?
\\
$(P)$ If $G$ is a cubic 3-connected graph, then $G$ has 
a $\Lambda $-packing that avoids at most two vertices of $G$, i.e. $\lambda (G) =   \lfloor v(G)/3 \rfloor $.}
\\[0.5ex]
\indent
We show, in particular, that claim $(P)$  
is equivalent to some seemingly 
stronger claims (see {\bf \ref{3-con}}).
For example, if $v(G) = 0 \bmod 3$, then the following claims are equivalent:
$G$ has a $\Lambda $-factor, 
for every  $e \in E(G)$ there exists a $\Lambda $-factor of $G$ avoiding (containing) $e$,
$G - X$ has 
a $\Lambda $-factor for every $X \subseteq E(G)$ such that  $|X| = 2$, and
$G - L$ has a $\Lambda $-factor for every
3-vertex path $L$ in $G$.
We also show that certain claims in theorem {\bf \ref{3-con}} are best possible. 
  We give a construction providing infinitely many cyclically 6-connected graphs $G$ with two disjoint 3-vertex paths L and L' such that $v(G) = 0 \bmod 3$ and 
  $G - L - L'$ has no $\Lambda $-factor.
It turns out  that if claim $(P)$ is true, then
Reed's dominating graph conjecture is true for cubic 
3-connected graphs.
\\[1ex]
 
\end{abstract}

\section{Introduction}

\indent

We consider undirected graphs with no loops and 
no parallel edges. All notions and facts on graphs, that are  
used but not described here, can be found in 
\cite{BM,D,Wst}.
\\[1ex]
\indent
Given graphs $G$ and $H$, 
an $H$-{\em packing} of $G$ is a subgraph of $G$ 
whose all components are isomorphic to $H$.
An $H$-{\em packing} $P$ of $G$ is called 
an $H$-{\em factor} if $V(P) = V(G)$. 
The $H$-{\em packing problem}, i.e., the problem of 
finding in $G$ an $H$-packing, having the maximum 
number of vertices, turns out to be $NP$-hard for every   
connected graph $H$ with at least three vertices \cite{HK}.
Let $\Lambda $ denote a 3-vertex path.
In particular, the $\Lambda $-packing problem
is $NP$-hard.
Moreover, the problem of packing $k$-vertex paths, $k \ge 3$, is $NP$-hard  for cubic graphs\cite{K1}.
Also, the 3-vertex path packing problem 
remains $NP$-hard even for cubic bipartite planar graphs
\cite{Kos3}.

Although the $\Lambda $-packing problem is $NP$-hard, 
i.e., possibly intractable in general, this problem turns out 
to be tractable for some natural classes of graphs. 
It would be  also interesting to find polynomial algorithms 
that would provide a good approximation  for 
the problem. Below 
(see {\bf \ref{eb(G)clfr}} -
{\bf \ref{2conclfr}}) there
are some examples of such results.
In each case the corresponding packing problem is 
polynomially solvable.
\\[1ex]
\indent
Let $v(G)$ and $\lambda (G)$ denote the number 
of vertices and the maximum number of disjoint
3-vertex paths in $G$, respectively.
Obviously, $\lambda (G) \le \lfloor v(G)/3 \rfloor $.

A graph is called {\em claw-free} if it contains no induced subgraph isomorphic
to $K_{1,3}$ (which is called a {\em claw}). 
A {\em block} of a graph $G$ is either an isolated vertex or a maximal connected subgraph $H$ of $G$ such that $H - v$ is connected for every vertex $v$ of $H$.
A block of a connected graph is called an 
{\em end-block} if it has at most  one vertex in common 
with any other block of the graph. 
Let $eb(G)$ denote the number of end-blocks of $G$.
\bs {\em \cite{KKN}}
\label{eb(G)clfr}
Suppose that $G$ is a connected claw-free graph and  
$eb(G) \ge 2$. 
Then  $\lambda(G) \ge \lfloor (v(G) - eb(G) + 2)/3 \rfloor$,
and this lower bound is sharp.
\es

\bs {\em \cite{KKN}}
\label{2conclfr} 
Suppose that $G$ is a connected and claw-free graph 
having at most two end-blocks 
$($in particular, a 2-connected and claw-free graph$)$.
Then $\lambda(G) = \lfloor v(G)/3 \rfloor$.
\es 

Obviously, the claim in {\bf \ref{2conclfr}} on claw-free graphs with exactly two end-blocks follows from  
{\bf \ref{eb(G)clfr}}.
\\[1ex]
\indent
In \cite{K,KM} we answered the following natural question:
\\[1ex]
\indent
{\em How many disjoint 3-vertex paths must a cubic graph have?}
\bs 
\label{km} If $G$ is a cubic graph, then 
$\lambda (G) \ge \lceil v(G)/4  \rceil$.
Moreover, there is a polynomial time algorithm for 
finding a $\Lambda $-packing having at least  
$\lceil v(G)/4  \rceil$ components.
\es

Obviously, if every component of $G$ is $K_4$, then
$\lambda (G) = v(G)/4$. Therefore the bound in {\bf \ref{km}} is sharp. There are some other connected cubic graphs for which the bound in {\bf \ref{km}} is attained.
\\[1ex]
\indent
Let ${\cal G}^3_2 $ denote the set of graphs with each vertex of degree at least $2$ and at most $3$.
In \cite{K} we  answered (among other results) the following  question:
\\[1ex]
\indent
{\em How many disjoint 3-vertex paths must an $n$-vertex  graph from ${\cal G}^3_2$ have?}
\bs 
\label{2,3-graphs} 
Suppose that $G \in {\cal G}^3_2$ and  
$G$ has no 5-vertex components.
Then $\lambda (G) \ge \lceil v(G)/4  \rceil $.
\es

Obviously, {\bf \ref{km}}  follows from {\bf \ref{2,3-graphs}} because if $G$ is a cubic graph, then $G \in {\cal G}^3_2$ 
and $G$ has no 5-vertex components. 
\\[1ex]
\indent
In \cite{K} we  also gave a construction  that allowed  to prove the following:
\bs 
\label{extrgraphs1}
There are infinitely many connected graphs for  which the bound  in {\bf \ref{2,3-graphs}} is attained.
Moreover, there are infinitely many subdivisions of
cubic 3-connected graphs for which the bound  in 
{\bf \ref{2,3-graphs}} is attained.
\es

The next interesting question is:
\\[1ex]
\indent
{\em How many disjoint 3-vertex paths must a cubic connected graph have?}
\\[1ex]
\indent
The following bound was recently found in \cite{Kos1}.
\bs 
\label{cubic-connected} 
If $G$ is a cubic connected graph and $v(G) \ge 17$, then 
$\lambda (G) \ge \frac{39}{152} v(G)$.
\es

The next natural question is:
\bs {\bf Problem.}
\label{2con-must}
How many disjoint 3-vertex paths must a cubic 2-connected graph have?
\es

\bs 
\label{cubic-connected} 
{\em \cite{Kos2}}
If $G$ is a cubic connected graph and $v(G) \ge 9$, then 
$\lambda (G) \ge \frac{3}{11} v(G)$.
\es

On the other hand, it is also natural to consider the following 
\bs {\bf Problem.}
\label{Pr2con}
Are there 2-connected cubic graphs $G$ such that
$\lambda (G) <   \lfloor v(G)/3 \rfloor $?
\es

In \cite{K2con-cbp}
we gave a construction that provided infinitely many   
2-connected, cubic, bipartite, and planar  graphs such that  
$\lambda (G) <   \lfloor v(G)/3 \rfloor $.
\\[1ex]
\indent
The main goal of  this paper (see also \cite{K3con-cub}) is to  discuss the following old open problem which is similar to Problem {\bf \ref{Pr2con}}.
\bs {\bf Problem.} {\em (A. Kelmans, 1984)}.
\label{Pr3con} 
Is the following claim true?
\\[1ex]
$(P)$ if $G$ is a 3-connected and cubic graph,
then $\lambda (G) =   \lfloor v(G)/3 \rfloor $.
\es

Originally we have
we've 
become interested in this problem because it is related with non-trivial relaxations of some graph Hamiltonicity problems we've been studying 
(see \cite{Kham}). 

Problem {\bf \ref{Pr3con}} is also related with 
 the minimum domination problem in a graph.
In \cite{R} B. Reed conjectured that
if $G$ is a connected cubic graph, then 
$\gamma (G) \le \lceil v(G)/3 \rceil $, where $\gamma (G)$ is the dominating number of $G$ (i.e., the size of a minimum vertex subset $X$ in $G$ such that every vertex in $G - X$ is adjacent to a vertex in $X$). It turns out that Reed's conjecture is not true for connected and even for 2-connected cubic graphs \cite{Kcntrex,KS}. 
If claim $(P)$ in {\bf \ref{Pr3con}} is true, then from 
{\bf \ref{3-con}} it follows, in particular, that Reed's conjecture is  true for 3-connected cubic graphs.
In \cite{KcubHam} we proved that if $G$ is a cubic Hamiltonian graph with $v(G) = 1\bmod 3$, then 
$\gamma (G) \le \lfloor v(G)/3 \rfloor $.

In this paper we show, in particular, that claim $(P)$ in 
{\bf \ref{Pr3con}} is equivalent to some seemingly 
stronger claims (see {\bf \ref{3-con}}).

In Section \ref{constructions} we give some notation, constructions, and simple observations.
In Section \ref{3connected} we formulate and prove our main theorem {\bf \ref{3-con}} concerning various claims  that are equivalent to claim $(P)$ in {\bf \ref{Pr3con}}.
We actually give different proofs of  {\bf \ref{3-con}}.
Thus, if there is  a counterexample $C$ to one of the above claims, then these different proofs 
provide different constructions 
of counterexamples to the other claims in {\bf \ref{3-con}}.
Additionally, 
different proofs provide better understanding of relations between various $\Lambda $-packing properties considered in {\bf \ref{3-con}}.
In Section \ref{AlmostCubic} we describe some 
properties of $\Lambda $-factors with respect to 
3-edge cut-matchings 
and triangles in cubic 3-connected graphs.
In Section \ref{homomorphism} we give 
two $\Lambda $-factor homomorphism theorems for cubic graphs.
Finally, in Section \ref{cycl6-connected}
we give a construction of cubic cyclically 6-connected graphs  showing that certain claims in 
{\bf \ref{3-con}} are best possible.

We can also prove that many equivalences in 
{\bf \ref{3-con}} are also true in the class of bipartite 3-connected cubic graphs and that if claim $(P)$ in {\bf \ref{Pr3con}} is true for bipartite graphs, then Reed's conjecture is also true for bipartite 3-connected cubic graphs.

\section{Notation, constructions, and simple observations}
\label{constructions}

\indent
  
We consider undirected graphs with no loops and 
no parallel edges unless stated otherwise. 
As usual, $V(G)$ and $E(G)$ denote the set of vertices and edges of $G$, respectively, and $v(G) = |V(G)|$.
If $X$ is a vertex subset or a subgraph of $G$, then let
$D(X,G)$ (or simply $D(X)$) denote the set of edges in 
$G$ having exactly one end-vertex in $X$, and let
$d(X,G) = |D(X,G)|$.
If $x \in V(G)$, then $D(x,G)$ is the set of edges in $G$ incident to $x$,
$d(x,G) = |D(x,G)|$, $N(x,G) = N(x)$ is the set of vertices 
in $G$ adjacent to $x$, and 
$\Delta (G) = \max \{d(x,G): x \in V(G)\}$.
If $e = xy \in E(G)$, then let $End(e) = \{x,y\}$.
Let $Cmp (G)$ denote the set of components of $G$ and $cmp(G) = |Cmp(G)|$.

\begin{figure}
  \centering
  \includegraphics{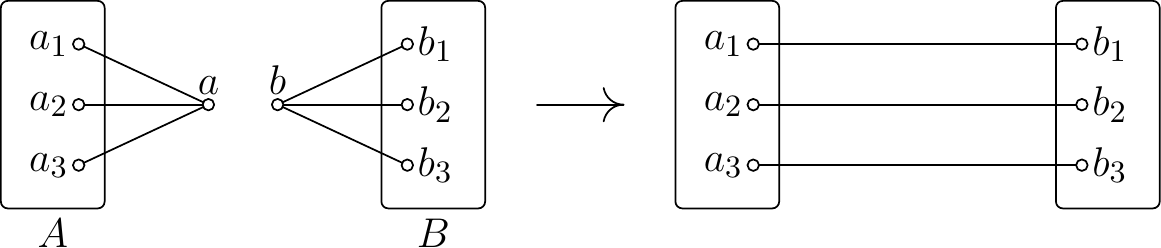}
  \caption{$Aa \sigma bB$}
  \label{fAasbB}
\end{figure}

We call an edge cut $K$ in a cubic graph a {\em blockade},
if $K$ is matching. A {\em $k$-blockade} is a blockade with $k$ edges.
Every connected cubic graph non-isomorphic to  
$K_4$ or $K_{3,3}$ has a blockade. 
For every such graph $G$, 
let $c(G)$ denote the size of a minimum blockade in 
$G$;
the parameter 
$c(G)$ is called  the {\em cyclic connectivity of} $G$ and
$G$ is called {\em cyclic $k$-connected}, if $c(G) \ge k$.
\\[1ex]
\indent
Let $A$ and $B$ be disjoint graphs, 
$a \in V(A)$,  $b \in V(B)$, and
$\sigma : N(a,A) \to N(b,B)$ be a bijection.
Let $Aa \sigma bB$ denote the graph
$(A - a) \cup (B - b) \cup \{x\sigma (x): x \in N(a,A)\}$.
We usually assume that 
$N(a,A) = \{a_1, a_2, a_3\}$, $N(b,B) = \{b_1, b_2, b_3\}$, and  $\sigma (a_i) = b_i$ for $i \in \{1,2,3\}$
(see Fig. \ref{fAasbB}). 
We also say that {\em $Aa \sigma bB$ is obtained from $B$ 
by replacing vertex $b$ by $(A - a)$ according to $\sigma $}.

Let $B$ be a cubic graph and $X \subseteq V(B)$. 
Let $A^v$, where $v \in X$, be a graph, 
$a^v$ be  a vertex of degree three in $A^v$, and 
$A_v = A^v - a^v$.
By using the above operation, we can build a graph
$G = B\{(A^v, a^v): v \in X\}$ by replacing each vertex $v$
of $B$ in $X$ by $A_v $ assuming that all $A^v$'s are disjoint. We call $G = B\{(A^v, a^v): v \in V(B)\}$ a {\em graph-composition},  $B$ -- the {\em frame}  and each 
$(A^v, a^v)$ -- a {\em brick} of this graph-composition. Let $D^v = D(A_v,G)$. If all $(A^v, a^v)$'s are isomorphic to  some
$(Z,z)$, then  instead of $B\{(A^v, a^v): v \in V(B)\}$ we simply write $B\{Z,z\}$.
For each $u \in V(B) \setminus X$, let $A^u$ be the graph having exactly two vertices $u$, $a^u$ and exactly three parallel edges connecting $u$ and $a^u$.
Then $G = B\{(A^v, a^v): v \in X\} = 
B\{(A^v, a^v): v \in V(B)\}$.
If, in particular, $X = V(B)$ and each $A^v$ is a copy of $K_4$, then $G$ is obtained from $B$ by replacing each vertex by a triangle.

Let $E' = E'(G) = E(G) \setminus \cup \{E(A_v):  v \in V(B)\}$.
Obviously, there is a unique bijection
$\alpha : E(B) \to E'(G)$ such that if $uv \in E(B)$, then 
$\alpha (uv)$ is an edge in $G$ having one end-vertex in
$A_u$ and the other in $A_v$.

Let $P$ be a $\Lambda $-packing in $G$.
For $uv \in E(B)$, $u \ne v$, we write  $u \neg ^p v$ or simply,
$u \neg  v$,
if $P$ has a 3-vertex path $L$ such that $\alpha (uv) \in E(L)$ and $|V(A_u) \cap V(L)| = 1$.
Let $P^v$ be the union of components of $P$ that meet
(i.e., have an edge in) $D^v$.
\\[1ex]
\indent
Obviously:
\bs
\label{AasbB} 
Let $k$ be an integer and $k \le 3$.
If $A$ and $B$ above are $k$-connected, cubic, 
bipartite, and planar graphs, then $Aa\sigma bB$ 
is also a $k$-connected, cubic, bipartite, and 
planar graph, respectively.
\es

From {\bf \ref{AasbB}} we have: 
\bs
\label{A(B)} 
Let $k$ be an integer and $k \le 3$.
If $B$  and each $A^v$ is a $k$-connected, cubic, 
bipartite, and planar graphs, then
$B\{(A^v, a^v): v \in V(B)\}$ is also a $k$-connected, 
cubic, bipartite, and planar graph, respectively.
\es

\begin{figure}
  \centering
  \includegraphics{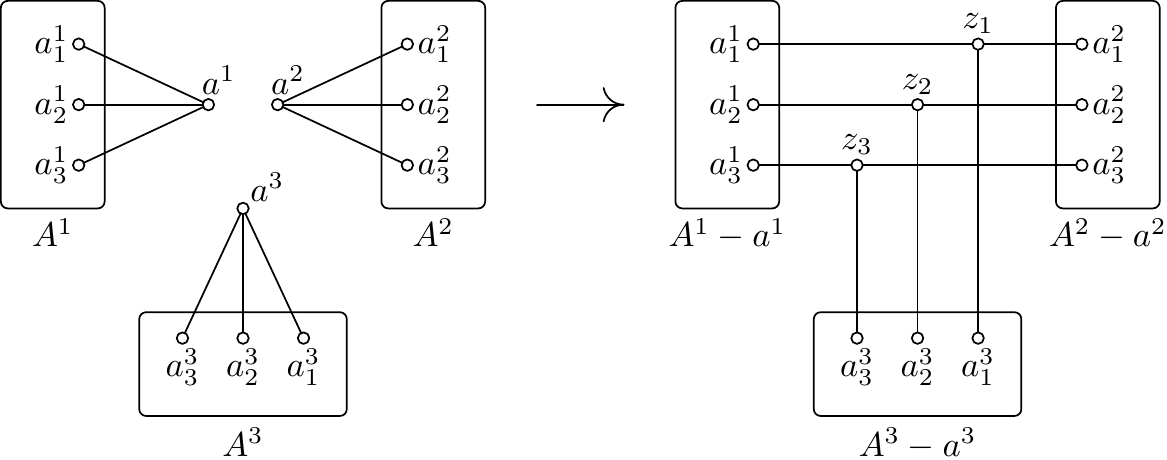}
  \caption{$Y(A^1,a^1; A^2,a^2; A^3,a^3)$}
  \label{fY}
\end{figure}

Let $A^1$, $A^2$, $A^3$ be three disjoint graphs,
$a^i \in V(A^i)$,  and $N(a^i,A^i) = \{a^i_1, a^i_2, a^i_3\}$, where $i \in \{1,2,3\}$.
Let  $F = Y(A^1,a^1; A^2,a^2; A^3,a^3)$  denote the graph obtained from 
$(A^1 - a^1) \cup (A^2 - a^2) \cup (A^3 - a^3)$ 
by adding three new vertices $z_1$, $z_2$, $z_3$ 
and the set of nine new edges $\{z_ja^i_j:  i, j \in \{1,2,3\}$
(see Fig. \ref{fY}). 
In other words, if $B = K_{3,3}$ is the complete
$(X,Z)$-bipartite graph with 
$X = \{x_1,x_2,x_3\}$ and  $Z = \{z_1,z_2,z_3\}$, then $F$ is obtained from  $B$ by replacing each vertex $x_i$ in $X$ by $A^i - a^i$ so that
$D(A^i - a^i, F) = \{a^i_jz_j: j \in \{1,2,3\}\}$.
Let $D^i = D(A^i- a^i, F)$.
If $P$ is a $\Lambda $-packing of $F$, then let 
$P^i = P^i(F)$ be the union of components of $P$ meeting $D^i$ and $E^i = E^i(P) = E(P) \cap D^i$, 
$i \in \{1,2,3\}$.

If each $(A^i, a^i)$ is a copy of the same $(A, a)$, 
then we write $Y(A, a)$ instead of 
$Y(A^1, a^1; A^2,a^2; A^3,a^3)$.
\\[1.5ex]
\indent
In particular, from {\bf \ref{A(B)} } we have: 
\bs
\label{Y(G1,G2,G3)} 
Let $k$ be an integer and $k \le 3$.
If each $A^i$ above is a $k$-connected, cubic, and 
bipartite graph, then $Y(A^1,a^1; A^2,a^2; A^3,a^3)$ 
{\em (see Fig. \ref{fY})} is also a $k$-connected, cubic, and bipartite graph, respectively.
\es

\begin{figure}
  \centering
  \includegraphics{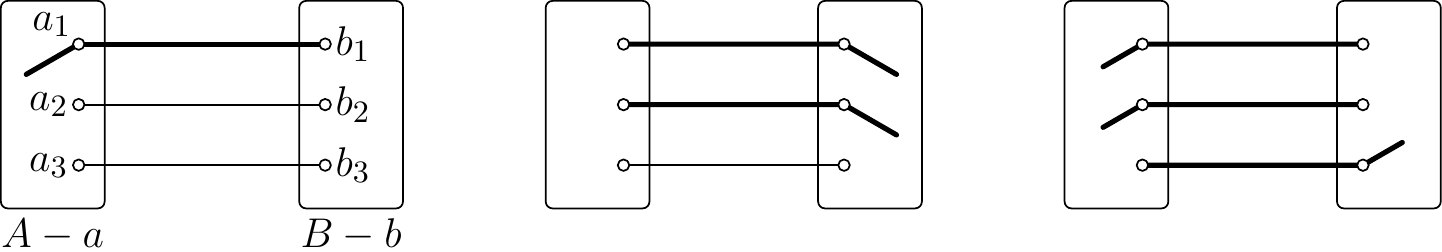}
  \caption{}
  \label{A1}
\end{figure}

We will also use the following simple observation.
\bs
\label{3cut}
Let $A$ and $B$ be disjoint graphs, 
$a \in V(A)$,  $N(a,A) = \{a_1, a_2, a_3\}$,  
$b \in V(B)$,  $N(b,B) = \{b_1, b_2, b_3\}$, and 
$G = Aa \sigma bB$, where each $\sigma (a_i) = b_i$
{\em (see Fig. \ref{fAasbB})}. 
Let $P$ be a $\Lambda $-factor of $G$ 
$($and so $v(G) = 0 \bmod 3$$)$ and 
$P'$ be the $\Lambda $-packing of $G$ consisting of the components $($3-vertex paths$)$ of $P$ that meet 
$\{a_1b_1, a_2b_2, a_3b_3\}$.
\\[1ex]
$(a1)$
Suppose that  $v(A) = 0 \bmod 3$, and 
so $v(B) = 2 \bmod 3$.  Then one of the following holds
{\em (see Fig. \ref{A1})}$:$
\\[0.5ex]
\indent
$(a1.1)$ $P'$ has exactly one component that
has two  vertices in $A - a$, that are adjacent 
$($and, accordingly, exactly one vertex in $B - b$$)$,
\\[0.5ex]
\indent
$(a1.2)$ $P'$ has exactly two components  and each component has exactly one vertex in $A - a$ 
$($and, accordingly, exactly two vertices in 
$B - b$, that are adjacent$)$,
\\[0.5ex]
\indent
$(a1.3)$ $P'$ has exactly three components $L_1$, $L_2$, $L_3$ and one of them, say $L_1$, has exactly one vertex in $A - a$ and each of the other two $L_2$ and $L_3$ has exactly two vertices in $A - a$, that are adjacent
$($and, accordingly, $L_1$ has exactly two  vertices in 
$B - b$, that are adjacent, and each of the other two $L_2$ and  $L_3$ has exactly  one vertex in $B - b$$)$. 
\\[1ex]
$(a2)$
Suppose that  $v(A) = 1 \bmod 3$, and 
so $v(B) = 1 \bmod 3$.  Then one of the following holds
{\em (see Fig. \ref{A2})}:
\begin{figure}
  \centering
  \includegraphics{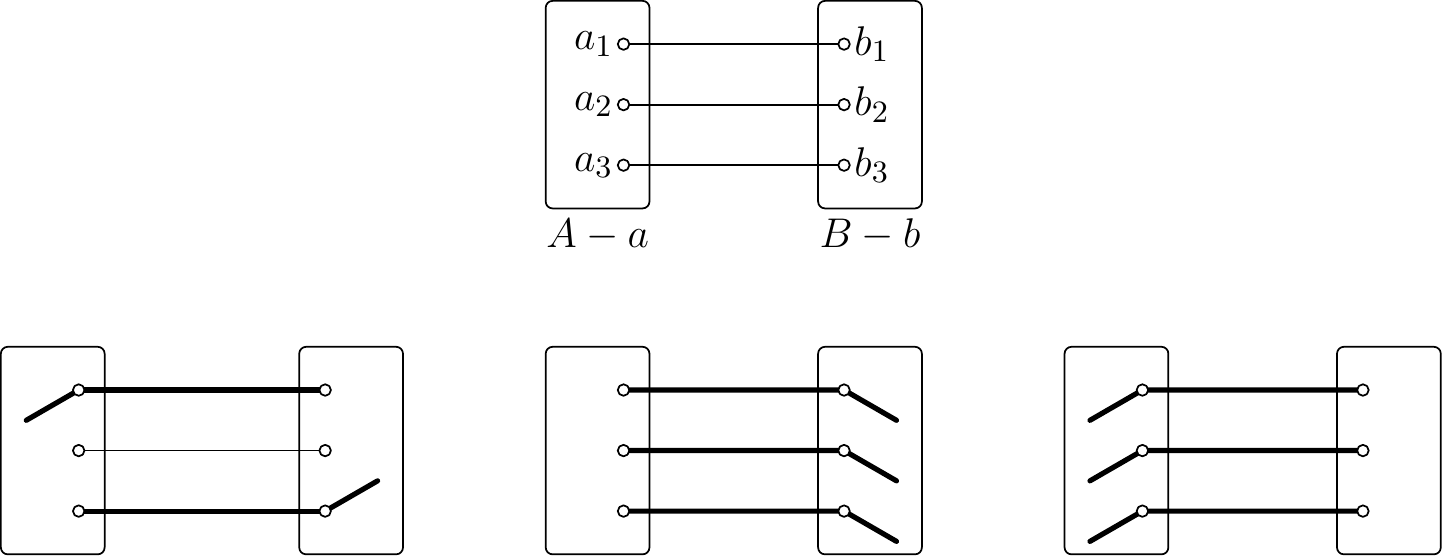}
  \caption{}
  \label{A2}
\end{figure}
\\[0.5ex]
\indent
$(a2.1)$ $P' = \emptyset $,
\\[0.5ex]
\indent
$(a2.2)$ $P'$ has exactly two components, say $L_1$ and $L_2$, such that  one of the them, say $L_1$, has exactly one vertex in 
$A - a$ and exactly two  vertices in $B - b$, that are adjacent, and the other component $L_2$ has exactly two vertices in 
$A - a$, that are adjacent, and exactly one vertex in $B - b$,
\\[0.5ex]
\indent
$(a2.3)$ $P'$ has exactly three components $L_1$, $L_2$, $L_3$ and either each $L_i$ has exactly one vertex in 
$A - a$ 
$($and, accordingly, has exactly two vertices in 
$B - b$, that are adjacent$)$ or each $L_i$ has exactly two vertices in $A - a$, that are adjacent 
$($and, accordingly, has exactly one vertex in $B - b$$)$.
\es

\section{$\Lambda $-packings in cubic 3-connected graphs}
\label{3connected}

\indent

The main goal of this section is to prove the following theorem showing that claim $(P)$ in 
{\bf \ref{Pr3con}} is equivalent to various seemingly 
stronger claims.
\bs
\label{3-con}
The following claims are equivalent for cubic 3-connected graphs $G$:
\\[1ex]
${\bf (z1)}$ 
$v(G) = 0 \bmod 6$ $\Rightarrow$ $G$ has 
a $\Lambda $-factor,
\\[1ex] 
${\bf (z2)}$ 
$v(G) = 0 \bmod 6$ $\Rightarrow$ for every 
$e \in E(G)$ there is a $\Lambda $-factor of $G$ 
avoiding $e$ $($i.e., $G - e$ has a $\Lambda $-factor$)$,
\\[1ex]
${\bf (z3)}$ 
$v(G) = 0 \bmod 6$ $\Rightarrow$ for every 
$e \in E(G)$ there exists a $\Lambda $-factor of $G$ 
containing $e$,
\\[1ex] 
${\bf (z4)}$
$v(G) = 0 \bmod 6$ $\Rightarrow$ for every 
$x \in V(G)$ there exists at least one 3-vertex path $L$ 
such that $L$ is centered at $x$ and $G - L$ has 
a $\Lambda $-factor, 
\\[1ex]
${\bf (z5)}$ 
$v(G) = 0 \bmod 6$ $\Rightarrow$ for every 
$x \in V(G)$ there  exist at least two 3-vertex paths $L$ 
such that $L$ is centered at $x$ and $G - L$ has 
a $\Lambda $-factor, 
\\[1ex]
${\bf (z6)}$ 
$v(G) = 0 \bmod 6$ $\Rightarrow$ for every 
$xy \in E(G)$ there  exist edges $xx', yy' \in E(G)$ such that
$G - xyy'$ and $G - x'xy$ have $\Lambda $-factors,
\\[1ex]
${\bf (z7)}$
$v(G) = 0 \bmod 6$ $\Rightarrow$ $G - X$ has 
a $\Lambda $-factor for every $X \subseteq E(G)$ such that  $|X| = 2$,
\\[1ex]  
${\bf (z8)}$ 
$v(G) = 0 \bmod 6$ $\Rightarrow$ 
$G - L$ has a $\Lambda $-factor for every
3-vertex path $L$ in $G$,
\\[1ex]  
${\bf (z9)}$
$v(G) = 0 \bmod 6$ $\Rightarrow$ for every 3-edge cut $K$ of $G$ and $S \subset K$, $|S| = 2$, there  
exists a $\Lambda $-factor $P$ of $G$ such that 
$E(P) \cap K = S$,
\\[1ex]
${\bf (t1)}$ 
$v(G) = 2 \bmod 6$ $\Rightarrow$ for every  
$x \in V(G)$ there 
exists $xy \in E(G)$ such that $G - \{x,y\}$ 
has a $\Lambda $-factor,
\\[1ex]
${\bf (t2)}$ 
$v(G) = 2 \bmod 6$ $\Rightarrow$ $G - \{x,y\}$ 
has a $\Lambda $-factor for every $xy \in E(G)$,
\\[1ex]
${\bf (t3)}$
$v(G) = 2 \bmod 6$ $\Rightarrow$ for every 
$x \in V(G)$ there  
exists a 5-vertex path $W$ such that $x$ is 
the center vertex of $W$ and $G - W$ has 
a $\Lambda $-factor 
{\em (see also {\bf \ref{3cut}} $(a1.2)$ and Fig. \ref{A1})},
\\[1ex]
${\bf (t4)}$
$v(G) = 2 \bmod 6$ $\Rightarrow$ for every 
$x \in V(G)$ and $xy \in E(G)$ there  
exists a 5-vertex path $W$
such that $x$ is the center vertex of $W$,
$xy \not \in E(W)$,  and $G - W$ has a $\Lambda $-factor
{\em (see also {\bf \ref{3cut}} $(a1.2)$ and Fig. \ref{A1})},
\\[1ex]
${\bf (f1)}$ 
$v(G) = 4 \bmod 6$ $\Rightarrow$ $G - x$ 
has a $\Lambda $-factor for every $x \in V(G)$,
\\[1ex]
${\bf (f2)}$ 
$v(G) = 4 \bmod 6$ $\Rightarrow$ $G - \{x, e\}$ 
has a $\Lambda $-factor for every $x \in V(G)$ and 
$e \in E(G)$,
\\[1ex]
${\bf (f3)}$
$v(G) = 4 \bmod 6$ $\Rightarrow$ for every 
$x \in V(G)$ there  
exists a 4-vertex path $Z$ such that
$x$ is an inner vertex of $Z$ and $G - Z$ has 
a $\Lambda $-factor
{\em (see also {\bf \ref{3cut}} $(a2.2)$ and Fig. \ref{A2})},
\\[1ex]
${\bf (f4)}$ 
$v(G) = 4 \bmod 6$ $\Rightarrow$ for every 
$x \in V(G)$ there 
exists $xy \in E(G)$ and a 4-vertex path $Z$ such that $x$ is an inner vertex of $Z$, $xy \not \in E(Z)$, 
and $G - Z$ has a $\Lambda $-factor 
{\em (see also {\bf \ref{3cut}} $(a2.2)$ and Fig. \ref{A2})},
\\[1ex]
${\bf (f5)}$
$v(G) = 4 \bmod 6$ $\Rightarrow$ for every 
$xy \in E(G)$ there exists a 4-vertex path $Z$ such that $xy$ is the middle edge of $Z$ 
and $G - Z$ has a $\Lambda $-factor 
{\em (see also {\bf \ref{3cut}} $(a2.2)$ and Fig. \ref{A2})},
\\[1ex]
${\bf (f6)}$
$v(G) = 4 \bmod 6$ $\Rightarrow$ for every $z \in V(G)$ and every 3-vertex path
$xyz$ there exists a 4-vertex path $Z$ such that $xyz \subset Z$, $z$ is an end-vertex of $Z$, 
and $G - Z$ has a $\Lambda $-factor 
{\em (see also {\bf \ref{3cut}} $(a2.2)$ and Fig. \ref{A2})}.   
\es

Theorem {\bf \ref{3-con}} follows from 
{\bf \ref{z1ottoz2}} -- {\bf \ref{z9ottoz1}} below.
In \cite{Kclfree} we have shown that claims $(z1)$ - $(z5)$ 
are true for cubic, 3-connected, and claw-free graphs.
In \cite{KcubHam} we have proved the following results related to $(t3)$, $(t4)$, and $(f3)$ - $(f6)$.

\bs
\label{Pi,W,Y}
Let $G$ be a cubic 2-connected graph. Then
\\[0.7ex]
$(a1)$ if $G - x$ has a $\Lambda $-factor for some vertex $x$ in $G$ {\em (and so $v(G) = 4 \bmod 6$)}, then
there exists a 4-vertex path $\Pi $ in $G$ such that
$G - \Pi $ has a $\Lambda $-factor,
\\[0.7ex]
$(a2)$ if $G - {x, y}$ has a $\Lambda $-factor for some edge $xy$ in $G$ {\em (and so $v(G) = 2 \bmod 6$)}, then
there exists a 5-vertex path $W$ and two disjoint 4-vertex paths $\Pi _1$ and $\Pi _2$ in $G$ such that
$G - W$ and $G - (\Pi _1 \cup \Pi _2)$ have 
$\Lambda $-factors,
and
\\[0.7ex]
$(a3)$ if $v(G) = 4 \bmod 6$ and $G$ has a Hamiltonian cycle, then there exists a claw $Y = K_{1,3}$ in $G$ such that $G - Y$ has a $\Lambda $-factor.
\es

From claim $(a3)$ in {\bf \ref{Pi,W,Y}} it follows that
if $G$ is a cubic Hamiltonian graph with 
$v(G) = 4 \bmod 6$, then its domination number is at most $\lfloor v(G)/3 \rfloor $.
\\[1ex]
\indent
The remarks below show that if claims $(z7)$, $(z8)$, 
$(t2)$ - $(t4)$, and 
$(f1)$ - $(f6)$
in {\bf \ref{3-con}} 
are true, then they are best possible in some sense.
\\[0.5ex]
\indent
$(r1)$ Obviously claim $(z7)$ is not true if 
condition ``$|X| = 2$'' is replaced by condition ``$|X| = 3$''.
Namely, if 
$G$ is a cubic 3-connected graph,
$v(G) = 0 \bmod 6$, 
$X$ is a 3-edge cut in $G$, and the two components of 
$G - X$ have the number of vertices not divisible by 3, 
then clearly $G - X$ has no  $\Lambda $-factor.
Moreover, for every 3-edge graph $Z$ with no isolated vertices and distinct from $K_{1,3}$, there exist infinitely many pairs $(G,X)$ such that $G$ is a cubic 3-connected graph, $X \subset E(G)$, $|X| = 3$, $G[X]$ (the subgraph of $G$ induced by $X$ is isomorphic to $Z$, and $G - X$ is connected and has no $\Lambda $-factor.
Also in Section \ref{AlmostCubic} 
we describe  an infinite set of cubic 3-connected graphs $G$  such that such that $G - E(T)$ has no  $\Lambda $-factor for every triangle $T$ in $G$.
\\[0.5ex]
\indent
$(r2)$ There exist infinitely many triples $(G, L, e)$ such that
$G$ is a cubic, 3-connected, bipartite, and planar graph,
$v(G) = 0 \bmod 6$, $L$ is a 3-vertex path in $G$,
$e \in E(G - L)$, and $(G - e) - L $ has no 
$\Lambda $-factor, and so claim $(z8)$ is tight.
Indeed, let  $v \in V(G)$,
$N(v,G) = \{x,y,z\}$, and $yz \in E(G)$, and so $vyz$ is a triangle. Since $G$ is 3-connected, $x$ is not adjacent to 
$\{y,z\}$. Let $L$ be a 3-vertex path in $G - v$
containing $yz$. Then $x$ is an isolated vertex in 
$(G - e) - L $, and therefore $(G - e) - L $ has no 
$\Lambda $-factor. Similarly, let $S = abcda$ be a 4-cycle in 
$G$, $e = xd$ be an edge in  a vertex in $G - E(S)$, and
$L = abc$ be  a 3-vertex path in $G$.
Then again $x$ is an isolated vertex in 
$(G - e) - L $, and therefore $(G - e) - L $ has no 
$\Lambda $-factor.  
 \\[0.5ex]
\indent
$(r2')$ The arguments in $(r2)$ also show that 
there exist infinitely many triples $(G, L, L')$ such that
$G - (L \cup L')$ has no $\Lambda $-factor, where
$L$ and $L'$ are disjoint 3-vertex paths in $G$,
and $G$ is a cubic, 3-connected, bipartite, and planar graph, and so claim $(z8)$ is tight in this sense as well.  Indeed, they can be obtained from the constructions in 
$(r2)$ by replacing vertex $x$ by a 3-vertex path in $G - L$ containing $x$.
 
In Section \ref{cycl6-connected} (see {\bf \ref{G,L,L'}}) we give a construction providing  
infinitely many cubic cyclically 6-connected graphs $G$ such that on the one hand,
$(G - e)- L$ has a $\Lambda $-factor for every 3-vertex path $L$ in $G$ and every edge in $G - L$, and on the other hand, $G$ has pairs of disjoint 3-vertex paths $L$ and  $L'$ such that $G - (L \cup L')$ has no $\Lambda $-factor. 
\\[0.5ex]
\indent
$(r3)$ 
There exist infinitely many triples $(G, xy, e)$ such that
$G$ is a cubic, 3-connected, bipartite, and planar 
graph, $v(G) = 2 \bmod 6$, $xy \in E(G)$, 
$e \in E(G - \{x,y\})$,
and $G - \{x, y ,e \}$ has no $\Lambda $-factor, and 
so claim $(t2)$ is tight.
\\[0.5ex]
\indent
$(r4)$ There exist infinitely many triples $(G, x, y)$ such that
$G$ is a cubic 3-connected graph, 
$v(G) = 2 \bmod 6$, $\{x,y\} \subset V(G)$, $x \ne y$, 
$xy \not \in E(G)$, and $G -  \{x,y\}$ has no 
$\Lambda $-factor, and so 
claim $(t2)$ is  not true if vertices $x$ and $y$ 
are not adjacent.
\\[0.5ex ]
\indent
$(r5)$ 
There exist infinitely many quadruples $(G,a,b, x)$ such that
$G$ is a cubic 3-connected graph with no 3-cycles and no 4-cycles, $v(G) = 4 \bmod 6$, $x \in V(G)$,
$a$ and $b$ are non-adjacent edges in $G - x$,
and $G - \{x, a,b \}$ has no $\Lambda $-factor, and 
so claim $(f2)$ is tight.
\\[0.5ex] 
\indent
$(r6)$ 
There exist  infinitely many triples $(G,L,x)$ such that $G$ 
is a cubic 3-connected graph 
with no 3-cycles and no 4-cycles, $v(G) = 4 \bmod 6$,
$x \in V(G)$, $L$ is a 3-vertex path in $G - x$, and 
$G - \{x, L \}$ has no $\Lambda $-factor (see claim $(f1)$).
\\[0.5ex] 
\indent
$(r7)$ 
There exist  infinitely many pairs $(G,\Pi)$ such that $G$ 
is a cubic 3-connected graph,
$v(G) = 4 \bmod 6$,
$\Pi$ is a 4-vertex path in $G$, and 
$G - \Pi$ has no $\Lambda $-factor 
(see claims 
$(f3)$ - $(f6)$).
\\[1ex]
\indent 
$(r8)$ 
There exist  infinitely many triples $(G,W)$ such that $G$ 
is a cubic 3-connected graph, 
$v(G) = 2 \bmod 6$,
$W$ is a 5-vertex path in $G$, and 
$G - W$ has no $\Lambda $-factor 
(see claims $(t3)$ and $(t4)$).
\\[1.5ex]
\indent 
The above discussion in $(r2)$ leads to the following question.

\bs {\bf Problem.}
Is the following claim true ?
Suppose that $G$ is a cubic 3-connected graph with no triangles and no squares $($or moreover, cyclically 5-connected$)$, $e$ is an edge in $G$, and $L$ is a 3-vertex path in $G - e$. Then $(G - e) - L$ has a $\Lambda $-factor.
\es

We need the following two results obtained before.
\bs {\em \cite{K2con-cbp}}
\label{Y,cmp(Pi)<3} 
Let $G = Y(A^1,a^1;  A^2,a^2; A^3,a^3)$ 
{\em (see Fig. \ref{fY})} and 
$P$  be a $\Lambda $-factor of $G$.
Suppose that each $A^i$ is a cubic graph and 
$v(A^i) = 0 \bmod 6$.
Then $cmp (P^i) \in \{1, 2\}$ for every $i \in \{1,2,3\}$.
\es
\noindent
{\bf Proof.}~
Let $i \in \{1,2,3\}$.
Since $D^i$ is a matching and $P^i$ consists of the components of $P$ meeting $D^i$, clearly $cmp (P^i) \le 3$.
Since $v(A^i) = -1 \bmod 6$, we have $cmp(P^i) \ge 1$.
It remains to show that $cmp (P^i) \le 2$.
Suppose, on the contrary, that $cmp(P^1) = 3$.
Since $P$ is a $\Lambda $-factor of $G$ and 
 $v(A^1 - a^1) = -1\bmod 6$, clearly
 $v(P^1) \cap V(A^1 - a^1) = 5$ and we can assume 
 (because of symmetry) that 
$P_1$ consists of three components  $a^1_3z_3a^2_3$,  
$z_1a^1_1 y^1$, and $z_2a^1_2 u^1$ for some 
$y^1, u^1 \in V(A^1)$ 
Then $cmp(P^3) = 0 $, a contradiction.
\ep
\bs {\em \cite{K2con-cbp}}
\label{Y,e-} 
Let $A$ be a graph, $e  = aa_1\in E(A)$, and 
$G = Y(A,a)$ {\em (see Fig. \ref{fY})}.
Suppose that
\\[0.5ex]
$(h1)$ $A$ is  cubic,
\\[0.5ex]
$(h2)$ $v(A) = 0\bmod 6$, and 
\\[0.5ex]
$(h3)$ $a$ has no $\Lambda $-factor containing 
$e = aa_1$.
\\[0.5ex]
\indent
Then $v(G) = 0 \bmod 6$ and $G$ has 
no $\Lambda $-factor.
\es

\bp (uses {\bf \ref{Y,cmp(Pi)<3}}).
Suppose, on the contrary, that $G$ has  a $\Lambda $-factor $P$. By definition of $G = Y(A,a)$, each $A^i$ is a copy of $A$ and edge $e^i = a^ia^i_1$ in $A^i$ is a copy of edge 
$e = aa_1$ in $A$ 
By {\bf \ref{Y,cmp(Pi)<3}}, $cmp(P^i) \in \{1, 2\}$
for every $i \in \{1,2,3\}$.
Since $P$ is a $\Lambda $-factor of $G$ and 
$v(A^i - a^i) = -1\bmod 6$, clearly
$E(P) \cap D^i $ is an edge subset of a $\Lambda $--factor 
of $A^i$ for every $i \in \{1,2,3\}$
(we assume that edge $z_ja^i_j$ in $G$ is edge $a^ia^i_j$ in $A^i$).
Since $a^1a^i_1$ does not belong to any
$\Lambda $-factor of 
$A^i$ for every $i \in \{1,2,3\}$, clearly
$E(P) \cap \{z_1a^1_1,z_1a^2_1,z_1a^3_1\}  = \emptyset $.
Therefore $z_1 \not \in V(P)$, and so $P$ is not 
a $\Lambda $-factor of $G$, a contradiction.
\ep 
\bs 
\label{z1ottoz2}
$(z1)$ $\Leftrightarrow$ $(z2)$.
\es
 
\bp (uses {\bf \ref{Y(G1,G2,G3)}} and {\bf \ref{Y,cmp(Pi)<3}}). 
Obviously, $(z1)$ $\Leftarrow$ $(z2)$. 
We prove $(z1)$ $\Rightarrow$ $(z2)$.

Suppose, on the contrary, that $(z1)$ is true but $(z2)$ 
is not true, i.e., there is a cubic 3-connected graph $A$ 
and $aa_1 \in E(G)$  such that $v(A) = 0 \bmod 6$ and
every $\Lambda $-factor of $G$ contains $aa_1$.
Let $G = Y(A^1,a^1; A^2, a^2; A^3, a^3)$, where 
each $(A^i, a^i)$ is a copy of $(A, a)$ and 
edge $a^ia^i_1$ in $A^i$ is a copy of edge $aa_1$ 
in $A$ 
(see Fig. \ref{fY}).
Since $A$ is cubic and 3-connected, 
by {\bf \ref{Y(G1,G2,G3)}},  $G$ is also cubic and 
3-connected.
Obviously, $v(G) = 0  \bmod 6$.
By $(z1)$, $G$ has a $\Lambda $-factor $P$.
Since each $v(A^i) = 0\bmod 6$, by {\bf \ref{Y,cmp(Pi)<3}},
$cmp(P^i)  \in \{1, 2\} $ for every $i \in \{1,2,3\}$.
Since $P$ is a $\Lambda $-factor of $G$ and 
$v(A^i - a^i) = -1\bmod 6$, clearly
$E(P) \cap D^i $ is an edge subset of 
a $\Lambda $-factor of $A^i$ for every $i \in \{1,2,3\}$
(we assume that edge $z_ja^i_j$ in $G$ is edge 
$a^ia^i_j$  in $A^i$).
Since $a^1a^i_1$ belongs to every  $\Lambda $-factor of 
$A^i$ for every $i \in \{1,2,3\}$, clearly
$ \{z_1a^1_1, z_1a^2_1,z_1a^3_1\} \subseteq E(P)$.
Therefore vertex $z_1$ has degree three in $P$, and so 
$P$ is not a $\Lambda $-factor of $G$, a contradiction.
\ep
\bs 
\label{z1ottoz3}
$(z1)$ $\Leftrightarrow$ $(z3)$.
\es
{\bf Proof.}~
Claim $(z1) \Leftarrow (z3)$ is obvious. 
Claim $(z1)\Rightarrow (z3)$ follows from {\bf \ref{Y,e-}}.
\ep
\bs 
\label{z1ottoz4}
$(z1)$ $\Leftrightarrow$ $(z4)$.
\es
{\bf Proof}~
(uses {\bf \ref{A(B)}} and {\bf \ref{3cut} }$(a1)$). 
Obviously, $(z1) \Leftarrow (z4)$.
We prove $(z1) \Rightarrow (z4)$.

\begin{figure}
  \centering
  \includegraphics{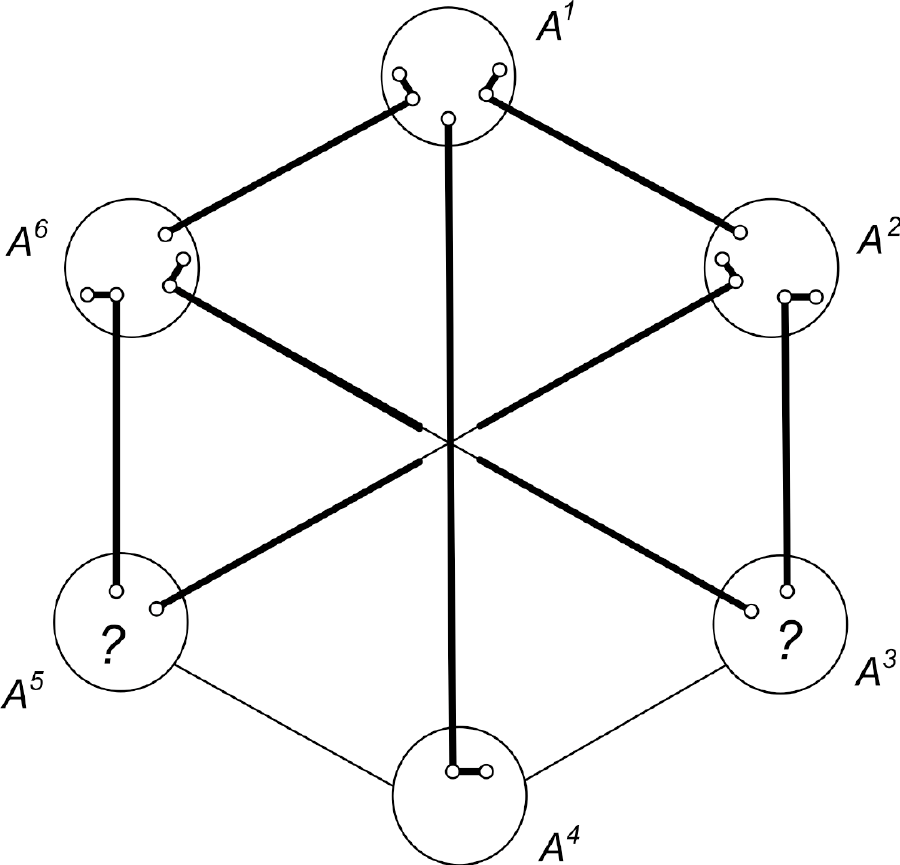}
  \caption{$(z1)$ $\Rightarrow$ $(z4)$}
  \label{fPrz1toz4}
\end{figure}

Suppose that $(z1)$ is true but  $(z4)$ is not true, 
i.e.,  there is a cubic 3-connected graph $A$
and $a\in V(A)$  such that $v(A) = 0 \bmod 6$ and
$a$ has degree one in every $\Lambda $-factor of $A$.
Let $G$ be the graph obtained from $B = K_{3,3}$ by replacing each vertex $v$ of $B$ by a copy $A_v$ of 
$A - a$ 
(see Fig. \ref{fPrz1toz4}). 
Obviously, $v(G) = 0 \bmod 6$ and by {\bf \ref{A(B)}},
$G$ is a cubic 3-connected graph.
By $(z1)$, $G$ has a $\Lambda $-factor $P$.
If $uv \in E(P)$, then let $L(uv)$ denote the component 
of $P$ containing $uv$. Let $V(B) = \{1, \ldots , 6\}$.

Since vertex $a$ has degree one in every $\Lambda $-factor of $A$ and each $(A^v, a^v)$ is a copy of $(A,a)$, by
{\bf \ref{3cut} }$(a1)$, we have $cmp(P^v) \in \{1, 3\}$.
\\[1ex]
${\bf (p1)}$
Suppose that there is $v \in V(B)$ such that 
$cmp(P^v) = 3$. By symmetry of $B$, we can assume that
$v = 1$ and,  by {\bf \ref{3cut}} $(a1.3)$,
$1 \neg 4$, $6 \neg 1$, and $2 \neg 1$
(see Fig. \ref{fPrz1toz4}). 
Let $x \in \{2,6\}$. Since $cmp(P^x) \in \{1,3\}$ and
$|V(L(1x) \cap V(A_x)| = 1$, clearly $cmp(P^x) = 3$ and,  
by {\bf \ref{3cut}} $(a1.3)$, $5 \neg x$, $3 \neg x$. 
Then by {\bf \ref{3cut} }$(a1.2)$, 
$cmp(P^s) = 2$ for $s \in \{3,5\}$, a contradiction.
\\[1ex]
${\bf (p2)}$ Now suppose that $cmp(P^v) = 1$ for every 
$v \in V(B)$. By symmetry, we can assume
$1 \neg 2$. Then $P^2$ contradicts {\bf \ref{3cut}} $(a1.1)$.
\ep
\bs
\label{z1ottot1}
$(z1)\Leftrightarrow (t1)$.
\es
{\bf Proof}~
(uses {\bf \ref{Y(G1,G2,G3)}}).
If $G = Y(A^1,a^1; A^2,a^2; A^3,a^3)$
(see Fig. \ref{fY}) and each $A^i$ is cubic 
and 3-connected, then
by {\bf \ref{Y(G1,G2,G3)}},  $G$ is also cubic and 
3-connected.
\\[1ex]
${\bf (p1)}$
We first prove $(z1) \Rightarrow (t1)$.
Suppose, on the contrary, that $(z1)$ is true but
$(t1)$ is not true, i.e., there is a cubic 3-connected graph 
$A$ and $a \in V(A)$ such that $v(A) = 2 \bmod 6$ and
 $A - \{a,y\}$ has no $\Lambda $-factor for every vertex 
 $y$ in $A$ adjacent to $a$. 
Let each $(A^i, a^i)$ above be a copy of $(A, a)$, and so
$A^i - \{a^i,a^i_j\}$ has no $\Lambda $-factor for every
$i,j \in \{1,2,3\}$. 
Obviously $v(G) = 0 \bmod 6$. By $(z1)$, $G$ has 
a $\Lambda $-factor $P$.
Then it is easy to see that  since each 
$v(A^i - a^i) = 1 \bmod 6$, 
there are $r,s, j \in \{1,2,3\}$ such that $r \ne s$ and 
$P^s = a^s_jz_ja^r_j$. 
Since $P$ is a $\Lambda $-factor of $G$, clearly 
$P \cap (A^s - \{z_j,a^s_j\})$ is a $\Lambda $-factor 
of $A^s - \{z_j,a^s\} = A^s - \{a^s,a^s_j\}$, a contradiction. 
\\[1ex]
${\bf (p2)}$
Now we prove $(z1) \Leftarrow (t1)$.
Suppose, on the contrary, that $(t1)$ is true but
$(z1)$ is not true, i.e., there is a cubic 3-connected
graph $A$ such that $v(A) = 0 \bmod 6$ and 
$A$ has no $\Lambda $-factor. 

Let $(A^i,a^i)$  be a copy of $(A,a)$ for 
$i \in \{1,2\}$, where $a \in V(A)$, and $(A^3,a^3)$ be 
a copy of $(H,h)$  for some cubic 3-connected graph 
$H$ and $h \in V(H)$, where $v(H) = 2 \bmod 6$. 
Obviously $v(G) = 2 \bmod 6$.
Suppose that $P$ is a $\Lambda $-factor of 
$G - \{z_1a^3_1\}$. Then $cmp(P^1) \le 2$. 
Since $v(A^1 - a^1) = -1 \bmod 6$, we have $cmp(P^1)\ge 1$.
Now since $P$ is a $\Lambda $--factor of 
$G- \{z_1a^3_1\}$ and $v(A^1- a^1) = -1\bmod 6$, 
clearly $E(P) \cap D^1 $ is an edge subset of 
a component of a $\Lambda $-factor of $A^1$
(we assume that edge $z_ja^1_j$ in $G$ is edge $a^1a^1_j$ in $A^1$).
Therefore $A$ has a $\Lambda $-factor, a contradiction.
\ep
\bs 
\label{z1ottof1}
$(z1)$ $\Leftrightarrow$ $(f1)$.
\es
{\bf Proof}~
(uses {\bf \ref{Y(G1,G2,G3)}}). 
By {\bf \ref{Y(G1,G2,G3)}}, if $G = Y(A^1,a^1; A^2,a^2; A^3,a^3)$ 
(see Fig. \ref{fY})
and each $A^i$ is cubic and 3-connected, then $G$ is also cubic and 3-connected.
\\[1ex]
${\bf (p1)}$
We first prove $(z1) \Rightarrow (f1)$.
Suppose, on the contrary, that $(z1)$ is true but $(f1)$ 
is not true, i.e., there is a cubic 3-connected graph $A$  
such that $v(A) = 4 \bmod 6$ and $A - a$ has no 
$\Lambda $-factor for some $a \in V(A)$.
Let each $(A^i,a^i)$ above be a copy of $(A, a)$.
Obviously, $v(G) = 0 \bmod 6$.
By $(z1)$, $G$ has a  $\Lambda $-factor $P$.
Let $P^i$ be the union of components in $P$ meeting $D^i$
and  $E^i = E^i(P)$.
 
Since $P$ is a  $\Lambda $-factor of $G$ and 
$v(A^i - a^i) = 3 \bmod 6$, clearly each $d(z_j, P) \le 2$ and 
each $|E^i| \in \{0, 2, 3\}$.
Since $A^i - a^i$ has no $\Lambda $-factor, 
each $|E^i| \in \{2, 3\}$, and so 
each $d(z_j, P) \ge 2$. Thus, each $d(z_j, P) = 2$, and so each $z_j$ is the center of a 3-vertex path in $P$ and each $|E^i| = 2$. Now, since $v(A^i - a^i) = 3 \bmod 6$, 
$P \cap (A^i - a^i - V(P^i))$ is not a $\Lambda $-factor of 
$A^i - a^i - V(P^i))$. Therefore $P$ is not a 
$\Lambda $-factor of $G$, a contradiction.
\\[1ex]
${\bf (p2)}$
Now we prove $(z1) \Leftarrow (f1)$.
Suppose, on the contrary, that $(f1)$ is true but $(z1)$ is 
not true, i.e., there is a cubic 3-connected graph $A$  
such that $v(A) = 0 \bmod 6$ and $A$ has no 
$\Lambda $-factor.
Let $(A^i,a^i)$ above be a copy of $(A, a)$ for
$i \in \{1,2\}$ and some $a \in V(A)$ and $(A^3,a^3)$ 
is a copy of $(H,h)$ for some cubic 3-connected graph 
$H$ and $h \in V(H)$, where $v(H) = 4 \bmod 6$ 
(see Fig. \ref{fY}).
Then $v(G) = 4 \bmod 6$. Let $x \in V(H - h)$.
Suppose that $G - x $  has a $\Lambda $-factor $P$.
Since $A$ has no $\Lambda $-factor, clearly 
$|E^1(P)| = |E^2(P)| = 3$. Then $P$ is not a 
$\Lambda $-factor of $G - x$, and so $(f1)$ is not true, 
a contradiction.
\ep
\bs 
\label{z4ottot2}
$(z4)$ $\Leftrightarrow$ $(t2)$.
\es
{\bf Proof}~
(uses {\bf \ref{Y(G1,G2,G3)}}, {\bf \ref{z1ottoz4}}, and 
{\bf \ref{z1ottot1}}).
Let $G = Y(A^1,a^1; A^2,a^2; A^3,a^3)$
(see Fig. \ref{fY}).
By {\bf \ref{Y(G1,G2,G3)}}, if each $A^i$ is 
cubic and 3-connected, then $G$ is also cubic and 
3-connected.
\\[1ex]
${\bf (p1)}$
We first prove $(z4) \Leftarrow (t2)$.
Obviously $(t2)  \Rightarrow (t1)$.
By  {\bf \ref{z1ottoz4}},
$(z1)$ $\Leftrightarrow$ $(z4)$.
By  {\bf \ref{z1ottot1}}, $(z1)\Leftrightarrow (t1)$.
The result follows.
\\[1ex]
${\bf (p2)}$ Now we prove
$(z4)\Rightarrow (t2)$.
Suppose, on the contrary, that $(z4)$ is true but $(t2)$ 
is not true. Then there exists a cubic 3-connected graph 
$A$ and $aa_1 \in E(G)$ such that $v(A) = 2 \bmod 6$ and 
$A - \{a,a_1\}$ has no $\Lambda $-factor.
Let each $(A^i, a^i)$ above be a copy of $(A,a)$ and 
edge $a^ia^i_1$ in $A^i$ be a copy of edge $aa_1$ in $G$.
Obviously $v(G) = 0 \bmod 6$.
Let $L^i = a^j_1z_1a^k_1$, where $\{i, j, k\} = \{1,2,3\}$.
By $(z4)$, $G$ has a $\Lambda $-factor $P$ containing 
$L^i$ for some $i \in \{1,2,3\}$, say for $i = 3$. 
If $s \in \{1,2\}$, then  $cmp (P^s) = 3$, because
$A^s - \{a^s, a^s_1\}$ has no $\Lambda $-factor and 
$v(A^s - a^s) = 1 \bmod 6$.
Also $cmp(P^3) \ge 1$, because
 $v(A^3 - a^3) = 1 \bmod 6$.
Then $P^1 \cup P^2 \cup P^3$ has at least four components each meeting $\{z_1,z_2,z_3\}$, 
a contradiction.
\ep 
\bs 
\label{z2toz5}
$(z2)$ $\Rightarrow$ $(z5)$.
\es
{\bf Proof}~
(uses  {\bf \ref{AasbB}} and {\bf \ref{3cut}} $(a1)$).
Suppose, on the contrary, that $(z2)$ is true but $(z5)$ 
is not true. Then there exists a cubic 3-connected graph 
$A$ and $a \in V(A)$ such that at most one 3-vertex path, centered at $a$  belongs to a  $\Lambda $-factor of $A$. It is sufficient to prove our claim in case when
$A$ has exactly one 3-vertex path, say $L = a_1 a a_2$, centered at $a$ and belonging to a  $\Lambda $-factor 
of $A$.  
Let $e_i = aa_i$, and so $E(L) = \{e_1,e_2\}$.

Let $B$ be the graph-skeleton of the three-prism: 
$V(B) = \{1, \ldots , 6\}$ and $B$ is obtained from two disjoint triangles $123$ and $456$ by adding three 
new edges $14$, $25$, and $36$.

\begin{figure}
  \centering
   \includegraphics{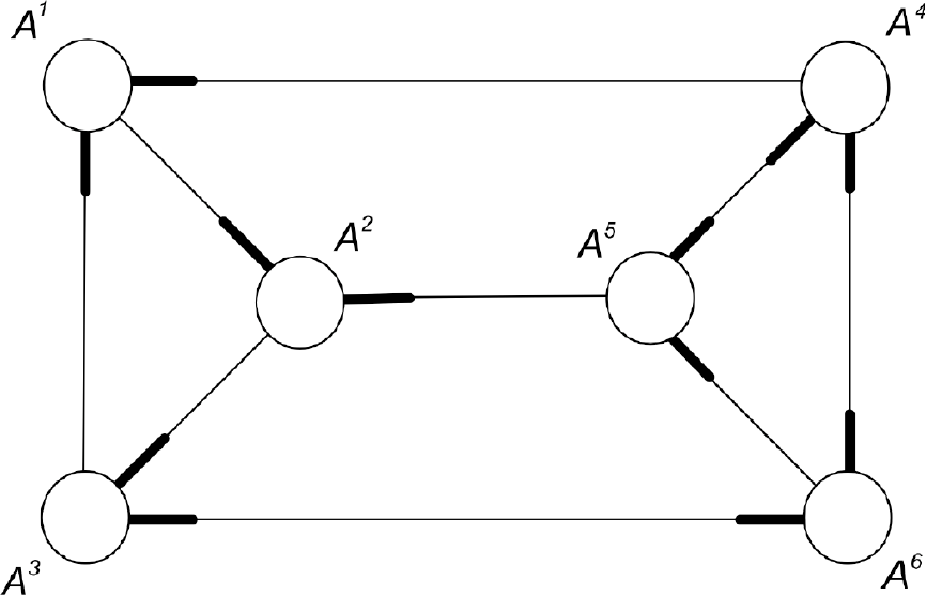}
  \caption{$(z2)$ $\Rightarrow$ $(z5)$}
  \label{fPrz2toz5}
\end{figure}

Let each $(A^v, a^v, a^v_1, a^v_2)$,  $v \in V(B)$ 
be a copy of $(A, a, a_1, a_2)$, 
and so edge $e^v_i = a^va^v_i$ in $A^v$ is a copy 
of edge $e_i = aa_i$ in $A$, $i \in \{1,2\}$. 
We also assume that all $A^v$'s are disjoint.
Let $G$ be a graph obtained from $B$ by replacing each 
$v \in V(B)$ by $A_v = A^v- a^v$
(see Fig. \ref{fPrz2toz5}).
We recall the following notation from Section \ref{constructions}:

\begin{itemize}
\item[(1)]
$E'(G) = E(G) \setminus \cup \{E(A_v): v \in V(B)\}$, 
\item[(2)]
$\alpha : E(B) \to E'(G)$ is the bijection such that if $uv \in E(B)$, then 
$\alpha (uv)$ is an edge in $G$ having one end-vertex in
$A_u$ and the other in $A_v$, and
\item[(3)]
if $P$ is  a $\Lambda $-packing in $G$, then 
for $uv \in E(B)$, $u \ne v$, we write  $u \neg ^p v$ (or simply, $u \neg  v$),
if $P$ has a 3-vertex path $L$ such that $\alpha (uv) \in E(L)$ and $|V(A_u) \cap V(L)| = 1$.
\end{itemize}

Given $v \in V(B)$, let  $S(v) $ be the set of two edges $e'_i$ in $E'(G)$ such that edge $e'_i$ is incident to vertex $a^v_i$ in $G$,  $i \in \{1,2\}$. 
We assume that each vertex $v$ in $B$ is replaced by 
$A_v$ (to obtain $G$) in such a way that
\\[0.5ex]
$S(1) = \{\alpha (13), \alpha (14\}$,
$S(2) = \{\alpha (21), \alpha (25)\}$,
$S(3) = \{\alpha (32), \alpha (36)\}$,
\\[0.5ex]
$S(4) = \{\alpha (45), \alpha (46)\}$,
$S(5) = \{\alpha (54), \alpha (56)\}$,
$S(6) = \{\alpha (63), \alpha (64)\}$.
\\[0.5ex]
In Figure \ref{fPrz2toz5} the edges in $S(v)$ are marked 
for every $v \in V(B)$.
\\[0.5ex]
\indent
By {\bf \ref{AasbB}}, $G$ is a cubic 3-connected graph. 
Since $v(B) = 0\bmod 6$, clearly also $v(G) = 0\bmod 6$.
By $(z2)$, $G' = G - \alpha (36)$ has a $\Lambda $-factor $P$.

We know that $A$ has exactly one 3-vertex path 
$L = a_1 a a_2$ centered at $a$ and belonging to 
a  $\Lambda $-factor of $A$ and that each 
$(A^v, a^v, a^v_1, a^v_2)$ is a copy of $(A, a, a_1,a_2)$, 
and so $v(A^v - a^v) = - 1 \bmod 6$.
Therefore, by {\bf \ref{3cut}} $(a1)$ and the assumption on $A$, 
the  $\Lambda $-factor $P$ satisfies the following condition
for every $v \in V(B)$: 
\\[1ex]
${\bf c(v)}$ if $cmp (P^v) = 2$, then 
$v  \neg a$ and $v  \neg b$, where 
$\{\alpha (va), \alpha (vb)\} = S(v)$.
\\[1ex]
\indent
Obviously, $|D^3| = |D^6| = 2$ in $G - \alpha (36)$.
Therefore $cmp (P^3) \le 2$ and $cmp (P^6) \le 2$.
Since $\alpha (36) \in S(3) \cap S(6)$, 
by conditions ${\bf c(3)}$ and ${\bf c(6)}$, 
$cmp (P^3) = cmp (P^6) = 1$.
Now by {\bf \ref{3cut}} $(a1)$, 
$x' \neg 3$ for some $x' \in \{1, 2\}$ and
$y' \neg 6$ for some $y' \in \{4, 5\}$.
\\[1ex]
${\bf (p1)}$ Suppose that $1\neg 3$.
Assume first that $\alpha (14) \not \in E(P)$.
Then $cmp(P^1) \le 2$. By {\bf \ref{3cut}} $(a1)$,
$cmp(P^1) = 2$. This contradicts ${\bf c(1)}$.
Thus we must have $\alpha (14) \in E(P)$.
\\[0.5ex]
${\bf (p1.1)}$ Suppose that $4 \neg 6$.

Suppose that $1\neg 4$. 
Then by {\bf \ref{3cut}} $(a1)$, $5 \neg 4$ and 
$5 \neg 2$. This contradicts ${\bf c(5)}$.

Now suppose that $4\neg 1$.
This contradicts ${\bf c(4)}$.
\\[0.5ex]
${\bf (p1.2)}$ Suppose that $5 \neg 6$.
Then $cmp (P^4) \le 2$.

Suppose that $1\neg 4$. 
By {\bf \ref{3cut}} $(a1.1)$, $cmp (P^4) = 1$.
Then $5 \neg 2$, which contradicts 
${\bf c(5)}$.

Now suppose that $4\neg 1$.
Then $cmp (P^4) = 2$, which contradicts ${\bf c(4)}$.
\\[1ex]
${\bf (p2)}$ Now suppose that $2\neg 3$.
Then $\alpha(13) \not \in E(P)$, and therefore $cmp(P^1) \le 2$.
By ${\bf c(2)}$ and {\bf \ref{3cut}} $(a1.3)$,
$1 \neg 2$ (and $5 \neg 2$).
Then by {\bf \ref{3cut}} $(a1.2)$, $cmp(P^1) = 2$,
which contradicts ${\bf c(1)}$.
\ep
\\

Since $(z7)$ $\Rightarrow$ $(z2)$, we obtain from 
{\bf \ref{z2toz5}}:

\bs 
\label{z7toz5}
$(z7)$ $\Rightarrow$ $(z5)$.
\es

\bs 
\label{z1ottoz6}
$(z1)$ $\Leftrightarrow$ $(z6)$.
\es

\noindent
{\bf Proof}.~
Obviously, $(z1)$ $\Leftarrow$ $(z6)$ and 
$(z5) \Rightarrow$ $(z6)$. 
As to $(z1)$ $\Rightarrow$ $(z6)$, it follows from 
{\bf \ref{z1ottoz2}} and {\bf \ref{z2toz5}}. 
\ep
\bs 
\label{z7toz8}
$(z7)$ $\Rightarrow$ $(z8)$.
\es

\noindent
{\bf Proof 1}.~
Suppose, on the contrary, that $(z7)$ is true but $(z8)$ 
is not true. Then there exists a cubic 3-connected graph $A$ and a 3-path $L = a_1aa_2$ in $A$ such that 
$v(A) = 0\bmod 6$ and $A - L$ has no  
$\Lambda $-factor. Let $N(a, A) = \{a_1, a_2, a_3\}$.
Let $(A^i;  a^i, a^i_1, a_2^i, a^i_3)$, $i \in \{1,2\}$, be two copies of $(A; a, a_1, a_2, a_3)$ and $A^1$, $A^2$ 
be disjoint graphs, and so $L^i = a^i_1 a^i a^i_2$ in 
$A^i$ is a copy of $L = a_1 a a_2$ in $A$.
Let $H = A^1a^1\sigma a^2A^2$, where 
$\sigma : N(a^1, A^1) \to N(a^2, A^2)$ is a bijection such that $\sigma (a^1_i) = a^2_i$ for $i \in \{1, 2, 3\}$.
Let $G$ be the graph obtained from $H$ by subdividing
edge $a^1_j a^2_j$ by a new vertex $v_j$ for $j \in \{1,3\}$ 
and by adding a new edge $v_1v_3$. 
Obviously, $G$ is a cubic 3-connected graph and
$v(G) = 0 \bmod 6$. 
By $(z7)$, $G' = G - \{a_3^1v_3, a_3^2 v_3\}$ has 
a $\Lambda $-factor $P$. Since $d(v_3,G') = 1$,  
$P$ has a 3-vertex path $L_1 = vv_1v_3$.
By symmetry of $G$, we can assume that 
$v = a_1^1$.  
Since $v(A^1 - a^1) = -1\bmod 6$,   $P$ has a 3-vertex path $L_2 = a_2^1a_2^2b$ for some $b \in V(A^2 - a^2)$.
Then $P \cap (A^1 - L)$ is a $\Lambda $-factor of 
$A^1 - L$, a contradiction.
\ep
\\[2ex]
{\bf Proof 2} (uses {\bf \ref{3cut}} $(a1)$).~
Suppose, on the contrary, that $(z7)$ is true but $(z8)$ 
is not true. Then there is a cubic 3-connected graph $A$
and a 3-path $L = a_1aa_2$ in $A$ such that 
$v(A) = 0\bmod 6$ and $A - L$ has no  $\Lambda $-factor. 
Let $N(a,A) = \{a_1,a_2,a_3\}$.

\begin{figure}
  \centering
   \includegraphics{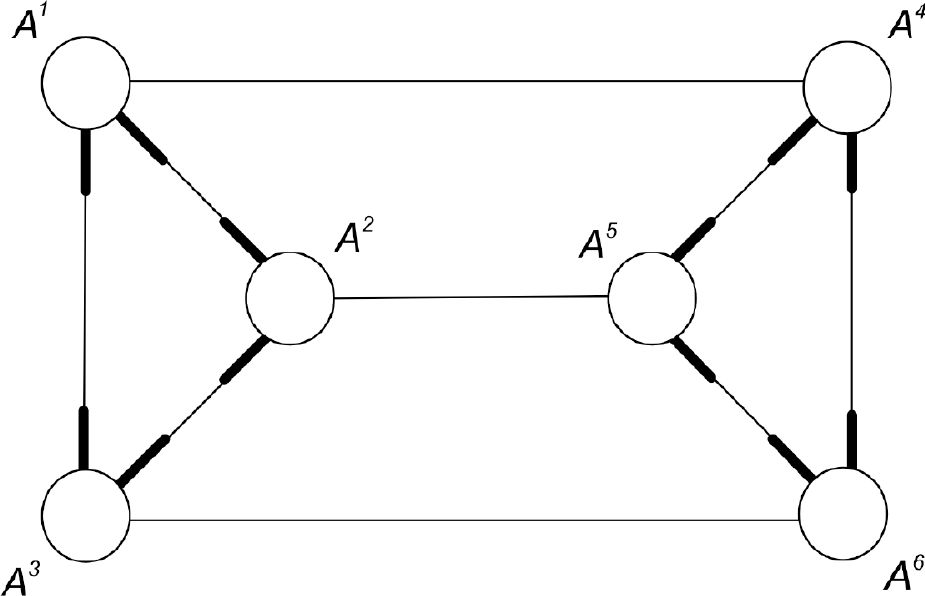}
  \caption{$(z7)$ $\Rightarrow$ $(z8)$}
  \label{fPrz7toz8}
\end{figure}

Let $B$,  $\{(A^v, a^v, a^v_1, a^v_2, a^v_3): v \in V(B)\}$, 
and $G$ be as in {\bf \ref{z2toz5}} (see Fig. \ref{fPrz7toz8}). 
Given $v \in V(B)$, let  $S(v) $ be the set of two edges 
$e'_i$ in $E'(G)$ such that edge $e'_i$ is incident to vertex $a^v_i$ in $G$,  $i \in \{1,2\}$. 
We assume that each vertex $v$ in $B$ is replaced by 
$A^v$ (to obtain $G$) in such a way that
\\[0.5ex]
$S(1) = \{\alpha (12), \alpha (13)\}$,
$S(2) = \{\alpha (21), \alpha (23)\}$,
$S(3) = \{\alpha (32), \alpha (31)\}$,
\\[0.5ex]
$S(4) = \{\alpha (45), \alpha (46)\}$,
$S(5) = \{\alpha (54), \alpha (56)\}$,
$S(6) = \{\alpha (64), \alpha (65)\}$.
\\[0.5ex]
In Figure \ref{fPrz7toz8}
the edges in $S(v)$ are marked 
for every $v \in V(B)$. 
\\[0.5ex]
\indent
Since $v(G) = 0\bmod 6$ and $G$ is cubic and 3-connected, by $(z7)$, $G - \{\alpha (14), \alpha (36)\}$ has a $\Lambda $-factor $P$. 
By {\bf \ref{3cut}} $(a1)$ and the lack of $\Lambda $-factor in each of $A^v - a_1^va^va_2^v$, 
the  $\Lambda $-factor $P$ satisfies the following condition  for every $v \in V(B)$: 
\\[1ex]
${\bf c(v)}$ if $cmp (P^v) = 2$, then 
$v  \neg a$ and $v  \neg b$, where 
$\{\alpha (va), \alpha (vb)\} \ne S(v)$.
\\[1ex]
\indent
Obviously, 
$|D^i| =  2$ in $G - \{\alpha (14), \alpha (36)\}$, 
and so  $cmp (P^i) \le 2$ for  $i \in \{1,3\}$.
Since $S(1) = \{12,13\}$ and $S(3) = \{31,32\}$, 
by conditions ${\bf c(1)}$ and ${\bf c(4)}$ we have 
$cmp (P^1) = cmp (P^3) = 1$.
Now by {\bf \ref{3cut}} $(a1.1)$, $2 \neg 1$ and $2 \neg 3$.
This contradicts ${\bf c(2)}$.
\ep

\bs 
\label{t2ottoz7}
$(t2)$ $\Leftrightarrow$ $(z7)$.
\es
{\bf Proof}~
(uses {\bf \ref{z1ottoz4}} and {\bf \ref{z4ottot2}}).
We first prove $(t2) \Rightarrow (z7)$.
Let $G$ be a cubic, 3-connected graph
with $v(G) = 0 \bmod 6$ and $a = a_1a_2$, $b = b_1b_2$ 
be two distinct edges of $G$.
Let $G'$ be the graph obtained from $G$ as follows:
subdivide edge $a_1a_2$ by a new vertex $a'$ and
edge $b_1b_2$ by a new vertex $b'$ and add a new edge 
$e = a'b'$. Then $G'$ is a cubic and 3-connected graph,
$v(G') = 2 \bmod 6$, and $G - \{a,b\} = G' - \{a',b'\}$.
By $(t2)$, $G' - \{a',b'\}$ has a $\Lambda $-factor.

Now we prove $(t2)\Leftarrow(z7)$.
Obviously, $(z7) \Rightarrow (z1)$.
By {\bf \ref{z1ottoz4}},
$(z1) \Rightarrow (z4)$ and by {\bf \ref{z4ottot2}},
$(z4)\Rightarrow (t2)$.
Implication $(t2)\Leftarrow(z7)$ also follows from obvious
$(z8)\Rightarrow (z4)$, from $(z7)\Rightarrow (z8)$, 
({\bf \ref{z7toz8}}), and from $(z4)\Rightarrow (t2)$
({\bf \ref{z4ottot2}}).
\ep
\\[1.5ex]
\indent
Here is a direct proof of $(z7)\Rightarrow (t2)$.
\bs 
\label{z7tot2}
$(z7)\Rightarrow (t2)$.
\es

\noindent
{\bf Proof.}~ 
Let $G$ be a cubic, 3-connected graph, 
$v(G) = 2\bmod 6$, $xy \in E(G)$, 
$N(x,G) = \{x_1,x_2, y\}$, and $N(y,G) = \{y_1,y_2, x\}$.
Let $G_1 = G - \{x, y\} \cup E_1$,
$G_2 = G - \{x, y\} \cup E_2$, and
$G_3 = G - \{x, y\} \cup E_3$,
where $E_1 = \{x_1y_1,x_2y_2\}$,  
$E_2 = \{x_1y_2,x_2y_1\}$, and 
$E_3 = \{x_1x_2,y_1y_2\}$.
Obviously, each $G_i$ is a cubic graph and $v(G_i) = 0 \bmod 6$.
It is easy to see that since $G$ is 3-connected, 
there is $s \in \{1,2,3\}$ such that $G_s$ 
is 3--connected. Clearly $G_s - E_s= G - \{x, y\}$.
By $(z7)$, $G_s - E_s$ has a $\Lambda $-factor.
\ep

\bs 
\label{z8ottoz7}
$(z8)$ $\Leftrightarrow$  $(z7)$.
\es
\noindent
{\bf Proof.}~
Obviously, $(z8) \Rightarrow (z4)$.
By {\bf \ref{z4ottot2}}, $(z4)\Rightarrow (t2)$ and
by {\bf \ref{t2ottoz7}}, $(t2) \Rightarrow (z7)$.
Therefore $(z8)$ $\Rightarrow$ $(z7)$.
By {\bf \ref{z7toz8}}, $(z7)$ $\Rightarrow$ $(z8)$.
\ep

\bs 
\label{z1ottoz8}
$(z1)$ $\Leftrightarrow$ $(z8)$.
\es

{\bf Proof.}~
 Obviously, $(z8)$ $\Rightarrow$ $(z1)$.
By {\bf \ref{z1ottoz4}},  $(z1)$ $\Rightarrow$ $(z4)$.
By {\bf \ref{z4ottot2}}, $(z4)$ $\Rightarrow$ $(t2)$.
By {\bf \ref{t2ottoz7}}, $(t2)$ $\Rightarrow$ $(z7)$.
By  {\bf \ref{z7toz8}}, $(z7)$ $\Rightarrow$ $(z8)$.
Therefore $(z1)$ $\Rightarrow$ $(z8)$.
\ep
\bs 
\label{z8tof1}
$(z8)$ $\Rightarrow$ $(f1)$.
\es

\noindent
{\bf Proof.}~
Let $G$ be a cubic 3-connected graph, 
$v(G) = 4 \bmod 6$, $x \in V(G)$, and 
$N(x,G) = \{x_1,x_2, x_3\}$.
Let $G'$ be the graph obtained from $G$ by replacing $x$ 
by a triangle $T$ with $V(T) =  \{x'_1, x'_2, x'_3\}$ so that $x_ix'_i \in E(G')$, $i \in \{1,2,3\}$.
Since $v(G) = 4 \bmod 6$, clearly $v(G') = 0 \bmod 6$.
Consider the 3-vertex path  $L' = x'_1x'_2x'_3$ in $G'$. 
By $(z8)$, $G' - L'$ has a $\Lambda $-factor  $P'$. Obviously, $P' - L'$ is 
a $\Lambda $-factor of $G' - L'$ and $G - x = G' - L'$.
\ep

\bs 
\label{z8tof2}
$(z8)$ $\Rightarrow$ $(f2)$.
\es

\begin{figure}
  \centering
   \includegraphics{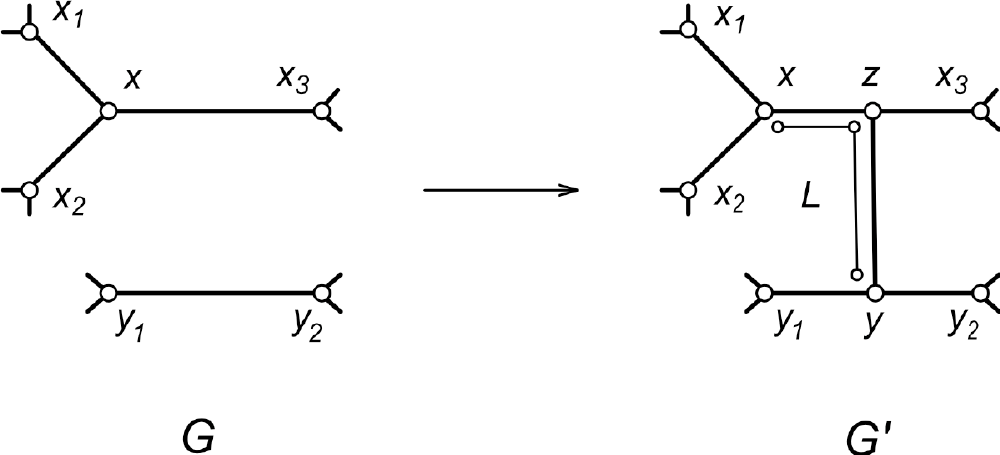}
  \caption{$(z8)$ $\Rightarrow$ $(f2)$}
  \label{fPrz8tof2}
\end{figure}

\bp (uses {\bf \ref{z8tof1}}).
Let $G$ be a cubic 3-connected graph, $v(G) = 4 \bmod 6$,
$x \in V(G)$, and $e = y_1y_2 \in E(G)$. We want to prove that if $(z8)$ is true, then $G - \{x,e\}$ has a 
$\Lambda $-factor. 
If $x  \in \{y_1,y_2\}$, then $G - \{x,e\} = G - x$, and therefore, by {\bf \ref{z8tof1}}, our claim is true. 
So we assume that $x  \not \in \{y_1,y_2\}$.
Let $N(x,G) = \{x_1,x_2,x_3\}$.
Let $G'$ be the graph obtained from $G$ by subdividing
edge $y_1y_2$ by a vertex $y$ and edge $xx_3$ by a vertex
$z$ and by adding a new edge $yz$
(see Fig. \ref{fPrz8tof2}).
Since  $v(G) = 4 \bmod 6$, clearly $v(G') = 0 \bmod 6$.
Since $x \not \in \{y_1,y_2\}$ and $G$ is cubic and 3-connected, $G'$ is also  cubic and 3-connected. Obviously,
$L = xzy$ is a 3-vertex path in $G'$ and $G - \{x,e\} =
G - \{x, y_1y_2\} = G' - L$. 
By $(z8)$, $G' - L$ has a $\Lambda $-factor.
\ep
\bs 
\label{f2tot4}
$(f2)$ $\Rightarrow$ $(t4)$.
\es

\noindent
{\bf Proof.}~
Let $G$ be a cubic, 3-connected graph, 
$v(G) = 2\bmod 6$, $x \in V(G)$, and 
$N(x,G) = \{x_1,x_2, x_3\}$.
Let $G'$ be the graph obtained from $G$ by replacing $x$ 
by a triangle $T$ with $V(T ) =  \{x'_1, x'_2, x'_3\}$ so that $x_ix'_i \in E(G')$, $i \in \{1,2,3\}$.
Since $v(G) = 2 \bmod 6$, clearly $v(G') = 4 \bmod 6$.
By $(f2)$, $G' - \{x'_i,x'_jx'_k\}$ has a $\Lambda $-factor
  $P_i$, where $\{i, j, k\} = \{1, 2, 3\}$. 
Since $x_jx'_j$ and $x_kx'_k$ are dangling edges in 
$G' - \{x'_i,x'_jx'_k\}$, clearly $x_jx'_j, x_kx'_k \in E(P_i)$ and 
$d(x'_j, P_i) = d(x'_k, P_i) =1$. Let $L_j$ and $L_k$ be the components of $P_i$ containing $x_jx'_j$ and $x_kx'_k$, respectively. Then $E(L_j) \cup E(L_k)$ induces in $G$ 
a 5-vertex path $W_i$ such that  $x$ is the center vertex of $W_i$ and $x_ix'_i \not \in E(W_i)$.
\ep
\bs 
\label{t3toz1}
$(t3)$ $\Rightarrow$ $(z1)$.
\es

\bp (uses {\bf \ref{3cut}} $(a1.3)$).
Suppose, on the contrary, that $(t3)$ is true but $(z1)$ 
is not true, i.e., there is a cubic 3-connected graph $A$  
such that $v(A) = 0 \bmod 6$ and $A$ has no 
$\Lambda $-factor. Let $a \in V(A)$.
Let $(A^i,a^i)$ above be a copy of $(A, a)$ for $i \in \{1,2\}$
and let $(A^3, a^3)$ be such that $v(A^3) = 2 \bmod 6$.
Let $G = Y(A^1,a^1; A^2,a^2; A^3,a^3)$
(see Fig. \ref{fY}).
Obviously, $v(G) = 2 \bmod 6$.
By $(t3)$, $G$ has a 5-vertex path $W$ such that $z_3$ is the center vertex of $W$ and $G - W$ has a 
$\Lambda $-factor $P$.
Obviously, $(A^i - a^i) \cap W = \emptyset $ for some 
$i \in \{1,2,3\}$.

Suppose that $(A^3 - a^3) \cap W = \emptyset $.
Then $W$ has an end-edge in $A^1 - a^1$ and 
in $A^2 - a^2$.
Since $A$ has no $\Lambda $-factor, 
by {\bf \ref{3cut}} $(a1.3)$, $D^i - e(W) \subseteq E(P)$ for 
$i \in \{1,2\}$. Then $P$ is not a $\Lambda $-factor of 
$G - W$, a contradiction.

Now suppose that $(A^3 - a^3) \cap W \ne \emptyset $.
By symmetry, we can assume 
that $(A^2 - a^2) \cap W = \emptyset $.
Then $W$ has an end-edge in $A^1 - a^1$ and 
in $A^3 - a^3$.
Then by {\bf \ref{3cut}} $(a1.3)$,  
$Cmp(P^1) = \{L_1,L_2\}$, where $L_1$  has an end-vertex 
in  $A^1 - a^1$ and $L_2$ has an end-edge in $A^1 - a^1$.
By symmetry, we can assume that $a^1_iz_i \in E(L_i)$ 
for $i \in \{1,2\}$. 
Then $L_1 = a^1_1z_1y$, where $y \in \{a^2_1,a^3_1\}$ and
$z_2$ is of degree one in $P$. 
Then, taking into account the values of $v((A^i - a^i)$ $i \in \{1,2,3\}$, $P$ is not a $\Lambda $-factor of 
$G - W$, a contradiction.
\ep
\bs 
\label{z8tof4}
$(z8)$ $\Rightarrow$ $(f4)$.
\es
\noindent
{\bf Proof.}~
Let $G$ be a cubic  3-connected graph, 
$v(G) = 4 \bmod 6$, $x \in V(G)$, and 
$N(x,G) = \{x_1,x_2, x_3\}$.
Let $G'$ be the graph obtained from $G$ by replacing $x$ 
by a triangle $ \Delta $ with 
$V(\Delta ) =  \{x'_1, x'_2, x'_3\}$, so that $x_ix'_i \in E(G')$, $i \in \{1,2,3\}$.
Since $v(G) = 4 \bmod 6$, clearly $v(G') = 0 \bmod 6$.
Since $G$ is cubic and 3-connected, $G'$ is also cubic and 3-connected.
Consider 3-vertex path  $L'_i = x_jx'_jx'_k$ in $G'$, 
where $\{i, j, k\} = \{1, 2, 3\}$. 
By $(z8)$, $G' - L'_i$ has a $\Lambda $-factor  $P'_i$.
Since $x_ix'_i$ is a dangling edge in 
$G' - L'_i$, clearly $x_ix'_i \in E(P'_i)$ and 
$d(x'_i, P_i) = 1$.
Let $\Lambda '_i$ be the 3-vertex path in $P'_i$ containing $x'_ix_i$, say $\Lambda '_i = x'_ix_iz$. 
Consider in $G$ the 4-vertex path 
$Z_k = x_kxx_iz$. Then $x$ is an inner vertex of $Z_k$, 
$xx_k \not \in E(Z_k)$, and $P'_i - \{L'_i \cup \Lambda '_i\}$ is a $\Lambda $-factor in $G - Z_k$.
\ep
\\[1.5ex]
\indent
Let $H'$ be a tree such that 
$V(H') = \{x,y\} \cup \{b^j: j \in \{1,2,3,4\}$ and 
$E(H') = \{xy, b^1x, b^2x, b^3y, b^4y\}$.
Let $H_i$, $i \in \{1,2,3\}$ be three disjoint copies of $H'$ 
with $V(H_i) = \{x_i,y_i\} \cup \{b^j_i: j \in \{1,2,3,4\}\}$.
Let $H$ be the graph obtained from these three copies by identifying 
for every $j$ three vertices $b^j_1$, $b^j_2$, $b^j_3$ with a new vertex $z^j$.
Let $A^i$, $i \in \{1,2,3,4\}$, be a cubic 
graph, $a^i \in V(A^i)$, and let
$G = H(A^1,a^1; A^2,a^2; A^3,a^3; A^4,a^4)$ be the graph obtained from $H$ by replacing each $z^j$ by $A^j - a^j$ assuming that all $A^i$'s are disjoint. 
\bs 
\label{f3toz1}
$(f3)$ $\Rightarrow$ $(z1)$. 
\es

\bp (uses {\bf \ref{A(B)}}).
Suppose, on the contrary, that $(f3)$ is true but $(z1)$ 
is not true, i.e., there exists a cubic 3-connected graph $A$
such that $v(A) = 0 \bmod 6$ and $A$ has no 
$\Lambda $-factor. Let $a \in V(A)$.
Let $A^4$ be any cubic 3-connected  graph with 
 $v(A^4) = 2 \bmod 6$ and each  
 $(A^i,a^i)$, $i \in \{1,2,3\}$, be a copy of $(A,a)$.
 
 Let the graph-composition $G = H(A^1,a^1; A^2,a^2; A^3,a^3; A^4,a^4)$ and its frame $H$ be as defined above.
Then $H$ is a cubic 3-connected graph. Since  each brick $A^i$ is  cubic and 3-connected, 
by {\bf \ref{A(B)}}, $G$ is also cubic and 3-connected.
Obviously, $v(G) = 4 \bmod 6$.
It is easy to see that $G - Z$ has no $\Lambda $-factor
for every 4-vertex path $Z$ in $G$ such that
$y_1$ is an inner vertex of $Z$. 
This contradicts $(f3)$.
\ep
\\[1.5ex]
\indent
Obviously, $(f4)$ $\Rightarrow$ $(f3)$. Therefore from 
{\bf \ref{f3toz1}} we have $(f4)$ $\Rightarrow$ $(z1)$.
Below we give a direct proof of this implication.
\bs 
\label{f4toz1}
$(f4)$ $\Rightarrow$ $(z1)$. 
\es
\noindent
{\bf Proof.}~
Suppose, on the contrary, that $(f4)$ is true but $(z1)$ 
is not true, i.e., there exists a cubic 3-connected graph $A$
such that $v(A) = 0 \bmod 6$ and $A$ has no 
$\Lambda $-factor. Let $a \in V(A)$.
Let $(A^3,a^3)$ be a copy of $(A, a)$
and let $(A^i, a^i)$ for $i \in \{1,2\}$ be copies of 
$(B,b)$, where $B$ is a  cubic 3-connected graph,
$v(B) = 2 \bmod 6$, and $b \in V(B)$.
Let $G = Y(A^1,a^1; A^2,a^2; A^3,a^3)$
(see Fig. \ref{fY}).
Obviously, $v(G) = 4 \bmod 6$.
By $(f4)$, $G$ has a 4-vertex path $Z$ such that $z_3$ 
is an inner  vertex of $Z$,  
$(A^3 - a^3) \cap Z = \emptyset $, and $G - Z$ has 
a $\Lambda $-factor, say $P$.
Since $A^3$ has no $\Lambda $-factor,
clearly $P$ is not a $\Lambda $-factor of $G - Z$, 
a contradiction.
\ep
\\[1.5ex]
\indent
Implication $(z8) \Rightarrow (f1)$ follows from obvious
$(z8) \Rightarrow$ $(z1)$ and from $(z1) \Rightarrow (f1)$,
by {\bf \ref{z1ottof1}}. 
It also follows from obvious
$(f2) \Rightarrow$ $(f1)$ and from $(z8) \Rightarrow (f2)$,
by {\bf \ref{z8tof2}}. Below we give a direct proof of this implication.
\bs 
\label{z1ottof6}
$(f6)$ $\Rightarrow$ $(f5)$ $\Rightarrow$ $(f4)$ $\Rightarrow$ $(z1)$ $\Rightarrow$ $(f6)$. 
\es

\bp (uses {\bf \ref{f4toz1}} and {\bf \ref{z1ottoz8}}).
Obviously, $(f6)$ $\Rightarrow$ $(f5)$ $\Rightarrow$ $(f4)$. By {\bf \ref{f4toz1}}, $(f4)$ $\Rightarrow$ $(z1)$. 
Therefore $(f6)$ $\Rightarrow$ $(z1)$. It remains to prove
$(z1)$ $\Rightarrow$ $(f6)$. 
By {\bf \ref{z1ottoz8}}, $(z1)$ $\Rightarrow$ $(z8)$.
Thus it is sufficient to show that $(z8)$ $\Rightarrow$ $(f6)$.
Let $G$ be a cubic 3-connected graph, 
$v(G) = 4 \bmod 6$, and $xyz$ is a 3-vertex path in $G$.
Let $N(y,G) = \{x,z,s\}$ and $G'$ be obtained from $G$ by
subdividing edges $yz$ and $ys$ by new vertices $z'$ and $s'$, respectively, and by adding a new edge $s'z'$.
Then $G'$ is a cubic 3-connected graph and $v(G') = 0 \bmod 6$. Consider the 3-vertex path  $L' = zz's'$ in $G'$. 
By  $(z8)$, $G'$ has a $\Lambda $-factor $P'$ containing $L'$. Since vertex $y$ has degree one in $G' - L'$, clearly
$P'$ has a 3-vertex path $Q' = yxq$. 
Let $Z$ be the 4-vertex path $qxyz$ in $G$.
Then $P' - (L' \cup Q')$ is a $\Lambda $-factor in $G - Z$.
\ep
\bs 
\label{z9ottoz1}
$(z9)$ $\Leftrightarrow$ $(z1)$.
\es

\noindent
{\bf Proof.}~
Obviously, $(z9)$ $\Rightarrow$ $(z8)$.
By the above claims,  $(z1)$, $(z8)$, $(t4)$, and $(f6)$ are equivalent. So we can use these claims to prove $(z9)$.
Let $G$ be a cubic 3-connected graph, $K$ a 3-edge cut of $G$, $S \subset K$ and $|S| = 2$, and $v(G) = 0 \bmod 6$. 
If the edges of $K$ are incident to the same vertex $x$ in $G$ (i.e., $K = D(x,G)$), then by $(z8)$, $G$ has  a $\Lambda $-factor $P$ of $G$ such that $E(P) \cap K = S$, and so our claim is true.
So we assume that the edges in $K$ are not incident to the same vertex in $G$. 
Then, since $G$ is cubic and 3-connected, clearly 3-edge cut  $K$ forms a matching. 
Let $A$ and $B$ be the two component of $G - K$. By the above arguments, we have 
$v(A)\ne 1$ and $v(B) \ne 1$.
Let $A^b$ be the graph obtained from $G$ by identifying the vertices of $B$ with a new vertex $b$ and similarly,
let $B^a$ be the graph obtained from  $G$ by identifying the vertices of $A$ with a new vertex $a$, and so $A = A^b - b$ and $B = B^a - a$. Then $D(b, A^b) = D(a,B^a) = D$, and
$S \subset D$. 
Let $S = \{a_1b_1, a_2b_2\}$, where 
$\{a_1,a_2\} \subset V(A)$ and $\{b_1,b_2\} \subset V(B)$.
Obviously,  $S$ forms a 3-vertex path $S_A = a_1ba_2$  in $A^b$ and a 3-vertex path $S_B = b_1ab_2$ in $B^a$. 
Since $G$ is cubic and 3-connected, both $A^b$ and $B^a$ are also cubic and 3-connected.
Since $v(G) = 0 \bmod 6$, there are two possibilities
(up to symmetry):
\\[0.5ex]
$(c1)$ $v(A^b) = 0 \bmod 6$ and $v(B^a) = 2 \bmod 6$ and 
\\[0.5ex]
$(c2)$ $v(A^b) = v(B^a) = 4 \bmod 6$.
\\[1ex]
\indent
Consider case $(c1)$. By $(z8)$, $A^b - S_A$ has a $\Lambda $-factor $P_A$. By $(t4)$, 
$B^a$ has a 5-vertex path $W$ such $a$ is the center vertex of $W$, $S_B \subset W$, and 
$B^a - W$ has a $\Lambda $-factor $P_B$.
Then $E(P_A) \cup E(P_B) \cup S$ induces 
a $\Lambda $-factor $P$ in $G$ such that 
$E(P) \cap K = S$.
\\[1ex]
\indent
Now consider case $(c2)$. 
By $(f6)$, we have:
\\[0.5ex]
$(a)$ 
$A^b$ has a 4-vertex path $Z_A$ such that $a_1$ is an end-vertex of $Z_A$, $S_A \subset Z_A$, and $A^b - Z_A $
has a $\Lambda $-factor, $P_A$;
\\[0.5ex]
$(b)$
likewise, $B^a$ has a 4-vertex path $Z_B$ such that $b_2$ is an end-vertex of $Z_B$, $S_B \subset Z_B$, and $B^a - Z_B$ has a $\Lambda $-factor $P_B$.
\\[0.7ex]
\indent
Then $E(P_A) \cup E(P_B) \cup S$ induces 
a $\Lambda $-factor $P$ in $G$ such that 
$E(P) \cap K = S$.
\ep

\section{Properties of $\Lambda $-factors in graph-compositions}
\label{AlmostCubic}

\indent

In this section we consider graph-compositions  having 
$\Lambda $-factors and with all bricks having the same number of vertices modulo 6.
It turns out that the $\Lambda $-factors of such 
graph-compositions intersect its 
3-blockades
in a special way and correspond to some special subgraphs in the  frame of the composition.

Let $G = B\{(A^v, a^v): v \in V(B)\}$ and $P$ be a $\Lambda $-factor of $G$ if any exists.
Let $F(P,G)$ be the subgraph in $B$ induced by the 
edge set $\alpha ^{-1}(E'(G) \cap E(P))$.
Let $\vec{F}(P,G)$ be the directed graph obtained from $F(P,G)$ by 
directing each edge $uv$ in $F(P,G)$ from $u$ to $v$, 
if  $u \neg ^p v$.
For a $\Lambda $-factor $P$ of $B$,
let
$\vec{P}$ denote the directed graph obtained from $P$ by directing the two edges of each 3-vertex path $L$  in $P$  to the central vertex of $L$. We call $\vec{P}$ a $\Lambda $-difactor of $B$.

A subgraph $C$ of a graph $H$ is called a {\em cycle-packing} in  $H$ if each component of $C$ is a cycle.
For a cycle-packing $Q$ in $B$,
 let
$\vec{Q}$ denote the directed graph obtained from $Q$  by directing its edges so that every cycle in $\vec{Q}$ becomes a directed cycle. We call $\vec{Q}$ a {\em cycle-dipacking} of $B$.
\bs 
\label{B(T)}
Let $B$ be a cubic 3-connected graph with no 
3-blockades.
Let $G$ be the graph obtained from $B$ by replacing each vertex by a triangle
{\em (i.e., $G = B\{K_4, x\}$, where $x \in V(K_4)$, 
and so $G$ is a claw-free graph).}
\\[0.7ex]
$(a1)$ Let  $P$ be a $\Lambda $-factor of $G$. Then
\\[0.5ex]
\indent
$(a1.1)$ $\vec{F}(P,G)$ is  a cycle-dipacking in $B$,
\\[0.5ex]
\indent
$(a1.2)$ 
$|K \cap E(P)| \in \{0, 2\}$ for every 
3-blockade $K$
of $G$, and 
\\[0.5ex]
\indent
$(a1.3)$ 
$E(T) \cap E(P) \ne \emptyset $ for every triangle $T$ of 
$G$.
\\[0.7ex]
$(a2)$ Let $\vec{Q}$ be a cycle-dipacking in $B$.
Then $G$ has a $\Lambda $-factor $P$ such that
$\vec{F}(P,G) = \vec{Q}$.
\es
{\bf Proof}~
(uses {\bf \ref{3cut}}).
Obiously, $(a1.1) \Rightarrow (a1.2)  \Rightarrow (a1.3)$.
So all we need prove $(a1)$.
Let $\vec{F} = \vec{F}(P,Q)$ and  $v \in V(\vec{F})$. 
Let $d_{in}(v)$ and $d_{out}(v)$ denote the number of edges coming in $v$ and going out of $v$ in $\vec{F}$, respectively.
By {\bf \ref{3cut}} and the fact that each vertex has been replaced by a triangle, there are three possibilities for every 
$v \in V(F)$:
\\[0.5ex]
$(c0)$ $d_{in}(v) = d_{out}(v) = 0$, 
\\[0.5ex]
$(c1)$ $d_{in}(v) = 1$ and $d_{out}(v) = 1$, and
\\[0.5ex]
$(c2)$ $d_{in}(v) = 0$ and $d_{out}(v) = 3$.
\\[0.5ex]
\indent
If $F$ has no vertices satisfying $(c2)$, then we are done.
So let $u$ be a vertex in $\vec{F}$ satisfying $(c2)$.
Let $S$ be a maximal directed path starting at $u$ and ending at $x$. Clearly, $d_{in}(x) = 1$. Then $x$ satisfies 
$(c1)$, and so $d_{out}(x) = 1$. Therefore $\vec{F}$ has an edge
$xy$ going out of $x$. Since $S$ is a maximal directed path
starting at $u$, clearly $y \in V(S - \{u,x\})$.
Then $d_{in}(y) = 2$, a contradiction.
Finally, $(a2)$ follows from the fact that  $G = B\{K_4, x\}$ and that every 2-edge matching covers the vertices of every triangle  in $K_4$.
\ep
\\[1ex]
\indent
Below we give a construction that provides a  big variety
of graph-compositions with the properties similar to that of 
a special graph-composition in {\bf \ref{B(T)}}. 
\bs
\label{v(A)=2mod6}
Let $B$ and each $A^u$, $u \in V(B)$, be a cubic 3-connected graphs and $G = B\{(A^u,a^u): u \in V(B)\}$. 
Suppose that 
\\[0.7ex]
$(h0)$ $v(B) =  0 \bmod 6$ and each $v(A^u) = 2 \bmod 6$, and so $v(G) =  0 \bmod 6$,
\\[0.7ex]
$(h1)$
each $A^u - (N(a^u) \cup a^u \cup x)$ has no $\Lambda $-factor
for every vertex $x$ in $A^u - (N(a^u) \cup a^u)$ adjacent 
to a vertex in $N(a^u)$,
and
\\[0.7ex]
$(h2)$ each $A^u - \{a^u, z\}$ has a $\Lambda $-factor for every  $a^uz \in E(A^u)$, and
\\[0.7ex]
$(h3)$
each $A^u - W$ has  a $\Lambda $-factor  for every 5-vertex path $W$ in $A^u$ centered at $a^u$.
\\[0.7ex]
\indent
Then
\\[0.7ex]
$(a1)$ if $P$ is a $\Lambda $-factor of $G$, then
$\alpha ^{-1}\{E(P) \cap E'(G)\}$ is the edge set of a 
$\Lambda $-factor $Q$ of $B$ and 
$\vec{F}(P,G) = \vec{Q}$,
\\[0.7ex]
$(a2)$ if $P$ is a $\Lambda $-factor of $G$ and if $B$ and each $A^u$ have no 
3-blockades,
then 
$|K \cap E(P)| \in \{1,2\}$ for every 
3-blockade
$K$ in $G$, and 
\\[0.7ex]
$(a3)$ if $\vec{Q}$ is a $\Lambda $-difactor  of $B$, 
then there exists a $\Lambda $-factor $P$ of $G$ such that 
$\vec{F}(P,G) = \vec{Q}$.
\es
{\bf Proof}~ (uses {\bf \ref{3cut}}).
Obviously, $(a2)$ follows from $(a1)$. We prove $(a1)$.
Let $E'(P) = E(P) \cap E'(G)$ and
$P^u$ be the union of the components (3-vertex paths) of $P$ meeting $D(A_u,G)$. By $(h0)$, $v(A^u) = 2 \bmod 6$. Therefore $P^u$ satisfies $(a1.1)$ or $(a1.2)$ or  $(a1.3)$ in {\bf \ref{3cut}}.
By $(h1)$, $P^u$ does not satisfy $(a1.3)$ in {\bf \ref{3cut}}.
Thus $\alpha ^{-1}(E'(P))$ is the edge set of a $\Lambda $-factor $P'$ in $B$ and the two edges of the same 3-vertex path $L$  in $P'$ are directed to the center vertex of $L$.
Finally, $(a3)$ follows from assumption $(h2)$ and $(h3)$.
\ep
\\[1ex]
\indent
It is also easy to prove the following claim similar to that above when the assumption $v(A^u) = 2 \bmod 6$ is replaced by $v(A^u) = 4 \bmod 6$. 
\bs
\label{v(A)=4mod6}
Let $B$ and each $A^u$, $u \in V(B)$, be a cubic 3-connected graphs and $G = B\{(A^u,a^u): u \in V(B)\}$. 
Suppose that 
\\[0.7ex]
$(h0)$ $v(B) =  0 \bmod 6$ and each $v(A^u) = 4 \bmod 6$, and so $v(G) =  0 \bmod 6$,
\\[0.7ex]
$(h1)$
each $A^u - T^u$  has no $\Lambda $-factor, where
$T^u$ is a tree subgraph of $A^u$ containing $a^u$ and such that $T^u - a^u$ is a 3-edge matching,
\\[0.7ex]
$(h2)$
each $A^u - a_u$ has a $\Lambda $-factor, and
\\[0.7ex]
$(h3)$ each $A^u - \Pi$ has a $\Lambda $-factor, where
$\Pi $ is a 4-vertex path in $A^u$ with $a^u$ being an inner vertex of  $\Pi $.
\\[0.7ex]
\indent
Then
\\[0.7ex]
$(a1)$ if $P$ is a  $\Lambda $-factor of $G$,  then 
$\vec{F}(P, G)$ is a cycle-dipacking in $B$,
\\[0.7ex]
$(a2)$ if $P$ is a  $\Lambda $-factor of $G$,  then
$|D(A_u,G) \cap E(P)| \in \{0, 2\}$ for every $u \in V(B)$, and
\\[0.7ex]
$(a3)$ if $\vec{Q}$ is a cycle-dipacking in $B$, then  $G$ has a  $\Lambda $-factor $P$ of $G$ such that 
$\vec{F}(P,G) = \vec{Q}$.
\es

The proof of $(a1)$ and $(a2)$ in {\bf \ref{v(A)=4mod6}} is similar to that in {\bf \ref{B(T)}}. The proof of  $(a3)$  in 
{\bf \ref{v(A)=4mod6}} is similar to that in {\bf \ref{v(A)=2mod6}}.
\\[1ex]
\indent
In the next two claims we give constructions of bricks satisfying the assumptions 
$(h0)$ - $(h3)$ in  {\bf \ref{v(A)=2mod6}} and 
{\bf \ref{v(A)=4mod6}}, respectively. 
\\[1ex]
\indent
Let $Z = Y(A^1,a^1; A^2,a^2; A^3,a^3)$
(see Fig. \ref{fY}), where $A^3$ is the graph having two vertices 
$z$,  $a^3$ and three parallel edges with the end-vertices
$z$,  $a^3$, and so $A^3 - a^3 = z$.
Let $K_i = D(A_i, Z)$, and so $K_i$ is the 
blockade
with three edges having exactly one end-vertex in  $A_i$, 
$i \in \{1,2\}$. 
\\[0.7ex]
\indent
Using {\bf \ref{3cut}}, it is easy to prove the following two claims.
\bs
\label{(Z,z)0} 
Let $(Z,z)$ be a pair defined above.
Suppose that $v(A^i) = 0 \bmod 6$ for $i \in \{1,2\}$.
Then $v(Z) = 2\bmod 6$ and   
\\[0.5ex]
$(a1)$
$Z - (N(z, Z)\cup z \cup v)$
has no $\Lambda $-factor for every vertex $v$ in 
$Z- (N(z, Z)\cup z)$ 
{\em $($see also {\bf \ref{3cut}} $(a1.3)$ and Fig. \ref{A1}$)$},
\\[0.5ex]
$(a2)$ if $zz' \in E(Z)$ and $P$ is a $\Lambda $-factor of 
$Z - \{z,z'\}$, then $|K_i \cap E(P)| \in \{1,2\}$ for 
$i \in \{1,2\}$, and
\\[0.5ex]
$(a3)$ if $W $ is a 5-vertex path in $Z$ with $z$ being the center vertex of  $W$ and $P$ is a $\Lambda $-factor of 
$Z - W$, then $|K_i \cap E(P)| \in \{1,2\}$ 
for $i \in \{1,2\}$.
\es

\bs
\label{(Z,z)4} 
Let $(Z,z)$ be a pair defined above.
Suppose that $v(A^i) = 4\bmod 6$ for  $i \in \{1,2\}$. Then $v(Z) = 4\bmod 6$ and  
\\[0.5ex]
$(a1)$  if $T$ is a tree subgraph of $Z$ containing $z$ and such that $T - z$ is a 3-edge matching, then $Z - T$ has no $\Lambda $-factor
{\em $($see also {\bf \ref{3cut}} $(a2.3)$ and Fig. \ref{A2}$)$},
\\[0.5ex]
$(a2)$ if $P$ is a $\Lambda $-factor of $Z - z$, then 
$|K_i \cap E(P)| \in \{0,2,3\}$ for $i \in \{1,2\}$, and if in addition, both $A^1$ and $A^2$ are isomorphic to $K_4$, then $|K_i \cap E(P)| = 2$ for $i \in \{1,2\}$,
\\[0.5ex]
$(a3)$ if $\Pi $ is a 4-vertex path in $Z$ with $z$ being an inner vertex of  $\Pi $ and $P$ is a $\Lambda $-factor of 
$Z - \Pi$, then $|K_i \cap E(P)| \in \{0,2\}$
for $i \in \{1,2\}$,
and
\\[0.5ex]
$(a4)$ if  $P$ is a $\Lambda $-factor of 
$Z - ( N(z)\cup z)$, then $|K_i \cap E(P)| = 0$ 
for $i \in \{1,2\}$.
\es

If $(Z,z)$ satisfies the assumptions in  
 {\bf \ref{(Z,z)0}}, then by  {\bf \ref{(Z,z)0}}, the graph-composition $B\{Z,z\}$ 
 satisfies the assumptions of  {\bf \ref{v(A)=2mod6}}, and therefore has  
the properties $(a1)$, $(a2)$, and $(a3)$ in 
{\bf \ref{v(A)=2mod6}}. 

Similarly, if $(Z,z)$ satisfies the assumptions in
 {\bf \ref{(Z,z)4}}, then by  {\bf \ref{(Z,z)4}}, the graph-composition $B\{Z,z\}$  satisfies the assumptions of  
 {\bf \ref{v(A)=4mod6}}, and therefore has  
the properties $(a1)$, $(a2)$, and $(a3)$ in 
{\bf \ref{v(A)=4mod6}}.  

  Notice that if any claim  in {\bf \ref{3-con}} is true, then there are infinitely many pairs $(Z,z)$ satisfying the assumptions  and, accordingly,
the conclusions of  {\bf \ref{(Z,z)0}}, as well as 
 of {\bf \ref{(Z,z)4}}.
 \\[1ex]
\indent
From {\bf \ref{v(A)=2mod6}}  and  {\bf \ref{v(A)=4mod6}} 
we have, in particular:
\bs 
\label{G,P,K} 
Let $\{i, j\}\in\{\{0,2\},\{1,2\}\}$.
Then there are infinitely many cubic 3-connected graphs 
$G$ 
such that $v(G) = 0 \bmod 6$ and
$|E(P) \cap K| \in \{i, j\}$
for every $\Lambda $-factor $P$ of $G$ and every  
3-edge cut $K$ of $G$.  
\es

\section{$\Lambda $-factor homomorphisms for graph-compositions}
\label{homomorphism}

\indent

Let, as in Section \ref{constructions}, 
$G = B\{(A^u,a^u): u \in V(B)\}$,
$A_u = A^u - a^u$,
and $E'(G) = E(G) \setminus \cup \{E(A_v):  v \in V(B)\}$.
Let $N(a^u, A^u) = N^u$.
As we mentioned in Section \ref{constructions}, there is a unique bijection
$\alpha : E(B) \to E'(G)$ such that if $uv \in E(B)$, then 
$\alpha (uv)$ is an edge in $G$ having one end-vertex in
$A_u$ and the other end-vertex in $A_v$.

Suppose that $P$ is a $\Lambda $-factor of $B$.
Let $V_s(P)$ be the set of vertices of degree $s$ in $P$ (and so $s \in \{1,2\}$).
For a 3-vertex path $uvw$ in $P$,
let
$A_1^u(P) = A^u - End(\alpha (uv))$ and
let $A_2^v(P)$ be the graph  obtained from $A_v$ by adding two edges  
$\alpha (uv)$ and $\alpha (vw)$ together with  their 
end-vertices.

Let $\Gamma (H)$ denote the set of $\Lambda $-factors of a graph $H$ and $\Gamma (G, P)$ denote the set of  
$\Lambda $-factors $Q$ of $G$ such that 
$\alpha ^{-1}(E'(G) \cap E(Q)\}) = E(P)$.
Let $X \bigotimes Y$ denote the Cartesian product of sets $X$ and $Y$ (we assume that if $Y = \emptyset $, then 
$X \bigotimes Y = X)$.
\\[1.5ex]
\indent
 From {\bf \ref{v(A)=2mod6}} we have:
\bs 
\label{Gamma}
Let $B$ and each $A^u$, 
$u \in V(B)$, be cubic 3-connected graphs.
Suppose that each $v(A^u)= 2 \bmod 6$ and 
each $(A^u, a^u)$ satisfies the following conditions: 
\\[0.7ex]
$(h1)$
$A^u - (N^u \cup a^u \cup y)$ has no $\Lambda $-factor
for every vertex $y$ in $A^u - (N^u \cup a^u)$ adjacent to a vertex in $N^u$
and
\\[0.7ex]
$(h2)$ $A^u - \{a^u, z\}$ has a $\Lambda $-factor for every
$a^uz \in E(A(u))$, 
and
\\[0.7ex]
$(h3)$ $A^u - W$
has a $\Lambda $-factor for every
5-vertex path $W$ in $A(u)$ centered at $a^u$.
\\[0.7ex]
\indent
Then $:$
\\[1ex]
$(\gamma 1)$
$\Gamma (G, P) \cap \Gamma (G, P') = \emptyset $ for 
$P, P' \in \Gamma (B)$, and $P \ne P'$,
\\[1ex]
$(\gamma 2)$
$\Gamma (G, P) = 
 [\bigotimes \{ \Gamma (A_1^v(P)): v \in V_1(P)\}] \bigotimes
[\bigotimes \{ \Gamma (A_2^u(P)): u \in V_2(P)\} ]$,
and
\\[1ex]
$(\gamma 3)$
$\Gamma (G) = \bigcup \{\Gamma (G, P): 
P \in \Gamma (B)\}$, 
\\[1ex] and so $G$ has a $\Lambda $-factor if and only if 
$B$ has a $\Lambda $-factor.
\es

Suppose that $\vec{C}$ is a cycle-dipacking in $B$.
If $uvw$ 
is a directed 3-vertex path in $\vec{C}$, then let
$A_*^v(\vec{C})$ be the graph obtained from $A_v$ by adding two edges  
$\alpha (uv)$ and $\alpha (vw)$ together with  their 
end-vertices, a new vertex $z$, and a new edge $wz$.

Let $\Gamma (G, \vec{C})$ denote the set of  $\Lambda $-factors $Q$ in $G$ such that  $\vec{F}(G,Q) = \vec{C}$.
Let ${\cal C}(B)$ denote the set of cycle-dipackings in $B$.
\\[1ex]
\indent
From {\bf \ref{v(A)=4mod6}} we have:
\bs 
\label{Gamma2}
Let $B$ and each $A^u$, 
$u \in V(B)$, be cubic 3-connected graphs.
Suppose that each $v(A(^u)= 2 \bmod 6$ and 
each $(A^u, a^u)$ satisfies the following conditions: 
\\[0.7ex]
$(h0)$ each $v(A^u) = 4 \bmod 6$, and so $v(B) =  0 \bmod 6$ and $v(G) =  0 \bmod 6$,
\\[0.7ex]
$(h1)$
each $A^u - T^u$  has no $\Lambda $-factor, where
$T^u$ is a tree subgraph of $A_u$ containing $a^u$ and such that $T^u - a^u$ is a 3-edge matching,
\\[0.7ex]
$(h2)$
each $A^u - a_u$ has a $\Lambda $-factor, and
\\[0.7ex]
$(h3)$ each $A^u - \Pi$ has a $\Lambda $-factor, where
$\Pi $ is a 4-vertex path in $A^u$ with $a^u$ being an inner vertex of  $\Pi $.
\\[0.7ex]
\indent
Then 
\\[1ex]
$(\gamma 1)$
$\Gamma (G, \vec{C}) \cap \Gamma (G, \vec{C}') = \emptyset $ for 
$\vec{C}, \vec{C}' \in {\cal C} (B)$, $\vec{C} \ne \vec{C}'$,
\\[1ex]
$(\gamma 2)$
$\Gamma (G, \vec{C}) = 
 [\bigotimes \{ \Gamma (A_*^v(\vec{C}): v \in V(\vec{C})\}]
 \bigotimes
[\bigotimes \{ \Gamma (A^u - a^u): u \in V(B) - V(\vec{C})\} ]$,
and
\\[1ex]
$(\gamma 3)$
$\Gamma (G) = \bigcup \{\Gamma (G,  \vec{C}): 
\vec{C} \in {\cal C} (B)\}$.
\es

As it is shown in Section \ref{AlmostCubic},
there are infinitely many pairs $(A,a)$ satisfying the assumptions in {\bf \ref{Gamma}}, as well as in {\bf \ref{Gamma2}}.

\section{Cyclically 6-connected cubic graphs with  special $\Lambda $-factor properties}

\label{cycl6-connected}

\indent

In this section we describe a sequence (mentioned in Section \ref{3connected}, remark $(r2)$) of cyclically 6-connected graphs $G$ with two disjoint 3-vertex paths $L$, $L'$ such that
$v(G) = 0 \bmod 6$ and $G - (L\cup L')$ has no $\Lambda $-factor.

We call a path $T$ in a graph $G$ a {\em thread} if the degree of each inner vertex of $T$ in $G$ is equal to two  and the degree of each end-vertex of $T$ in $G$ is not equal to two. 

\begin{figure}
  \centering
  \includegraphics{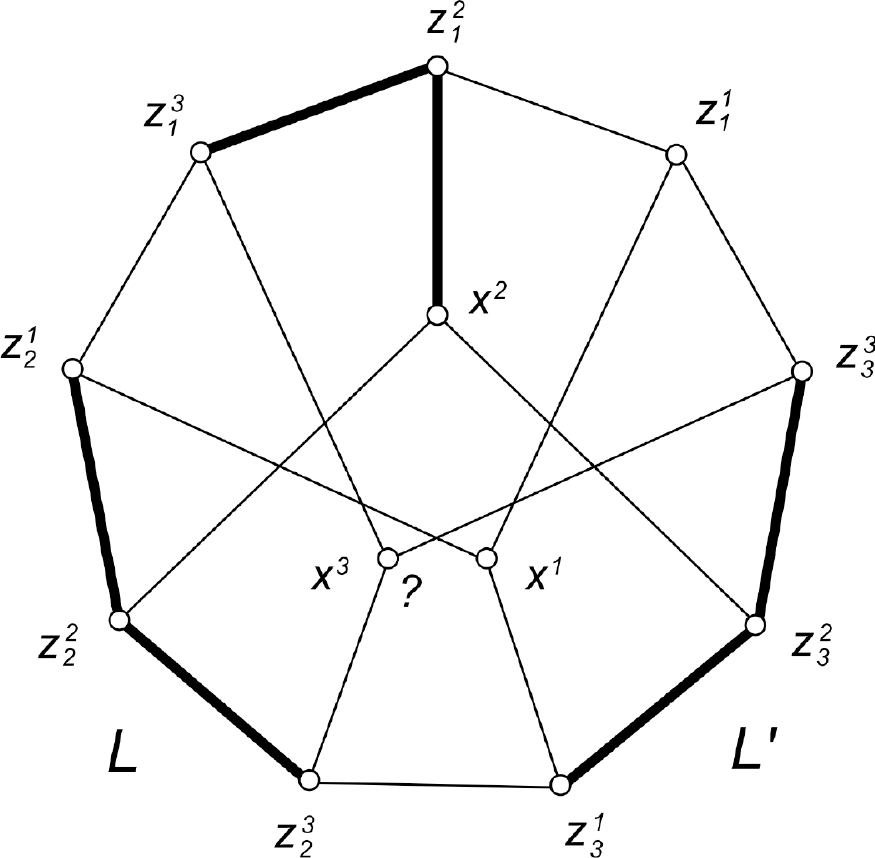}
    \caption{$(R_1, L, L')$, $v(R_1) = 12$}
  \label{fR1}
\end{figure}

Let $C_s$ be a cycle with $9s$ vertices, $s \ge 1$, and let
$\{L_k: k \in \{1, \ldots , 3s\}\}$ be a $\Lambda $-factor 
\\[0.1ex]
of $C_s$, where $L_i = (z_i^1z_i^2 z_i^3)$.
Let $R_s$ be the graph obtained from $C_s$ by adding the set 
\\[0.1ex]
$\{x^j_i:  i \in \{1, \ldots , s\}, j \in \{1,2,3\}\}$ 
of $3s$ new vertices and  the set 
$\{x^j_iz^j_i, x^j_iz^j_{i+s}, x^j_iz^j_{i+2s}:  
\\[0.1ex]
i \in \{1, \ldots , s\}, j \in \{1,2,3\}\}$ of $9s$ new edges 
(see Fig. \ref{fR1} for an illustration of $(R_1, L, L')$).
Let $X_i = \{x_i^1,x_i^2,x_i^3\}$.

Let  
$Q^i_{s-1} = R_s - X_i$.
Then the vertices of 3-vertex paths $L_i$, $L_{i+s}$, and $L_{i+2s}$ are of degree two in $Q^i_{s-1}$.
Let $T_i$, $T_{i+s}$, and $T_{i+2s}$ be the threads in $Q^i_{s-1}$ containing  $L_i$, $L_{i+s}$, and $L_{i+2s}$, respectively, and so each of these threads has exactly five vertices.
Let $R^i_{s-1}$ be the graph obtained from $Q^i_{s-1}$ by replacing threads $T_i$, $T_{i+s}$, and $T_{i+2s}$ by  edges $t_i$, $t_{i+s}$, and $t_{i+2s}$. Then obviously,
$R^i_{s-1}$ is isomorphic to $R_{s-1}$.  

Let $Y_i$ be the subgraph of $R_s$ induced by 
the vertex subset $V(L_i\cup L_{i+s}\cup L_{i+2s}) \cup 
 X_i$. 
 Then $R_s = Q^i_{s-1} \cup Y_i$, 
 $C_s \cap Y_i = L_i\cup L_{i+s}\cup L_{i+2s}$, and
 the set $K^i$ of edges having one vertex in $Y_i$ and the other vertex in $R_s - Y_i$ is a 
6-blockade in $R_s$, $s \ge 2$.
\bs
\label{G,L,L'} 
Let $R_s$, $s \ge 1$, be the graph defined above, and so 
$v(R_s) = 12 s$.
Then 
\\[0.5ex]
$(a1)$
$R_s - (L \cup L')$ has no $\Lambda $-factor for every pair
$\{L,L'\} \subset \{L_i,L_{i+s}, L_{i+2s}\}$, 
$i \in \{1, \ldots , s\}$,  
\\[0.7ex]
$(a2)$ $c(R_s) = 6$ for $s \ge 2$ and 
$c(R_1) = 5$ $($and for every 5-blockade
$M$ of $R_1$ one of the two components of $R_1 - M$ 
is a 5-cycle$)$, and  
\\[0.7ex]
$(a3)$ $R_s - L - e$ has a $\Lambda $-factor for every edge $e$ and every 3-vertex path $L$ in $R_s - e$. 
\es 

\noindent
{\bf Proof.}~ The proof will consist of three parts.
\\[1ex]
${\bf (p1)}$
We prove $(a1)$.
By symmetry, we can assume that 
$L = L_{1+s}$ and $L' = L_{1+2s}$. Let $G = R_s - (L \cup L')$. Suppose, on the contrary, that
$G = R_s - (L \cup L')$ has a $\Lambda $-factor $P$.
Then each edge $x_1^jz_1^j$, $ j \in \{1,2,3\}$, is a dangling
edge in $G$. Then  $x_1^2z_1^2$ belongs to a 3-vertex path  $L''$. By symmetry, we can assume that 
$L'' = z_1^3z_1^2x_1^2$. Then $x_1^1$ is an isolated vertex in $G - L''$ (see, for example, Fig. \ref{fR1}, where $x^i = x^i_1$). Therefore $P$ is not a $\Lambda $-factor of $G$, a contradiction.
\\[1ex]
${\bf (p2)}$
We prove $(a2)$ by induction on $s$. 
It is easy to check that $R_1$ has the properties described in $(a2)$. So let $s \ge 2$.
Since $R_s$ has 6-cycles, $c(R_s) \le 6$.
So our goal is to prove that  $c(R_s) = 6$.
Let $M$ be a minimum 
blockade
in $R_s$.
Obviously, each $R_s$ is a connected graph, and so 
$M \ne \emptyset $.

Let $Q_{s-1} = Q^s_{s-1}$, and so the vertices of 3-vertex paths $L^1 = L_s$, $L^2 = L_{2s}$, and $L^3 = L_{3s}$
are of degree two in $Q_{s-1}$. 
Let $T^1 = T_s$, $T^2 = T_{2s}$, and $T^3 = T_{3s}$ be the threads in $Q_{s-1}$ containing  
$L^1 $, $L^2  $, and $L^3 $, respectively.
Then graph $R_{s-1}$ can be obtained from $Q_{s-1}$ by replacing each thread $T^j$ by a new edge $t^j$, 
$j \in \{1,2,3\}$. 
Let $Y = Y_s$.
 Then $R_s = Q_{s-1} \cup Y$ and 
 $C_s \cap Y = L^1 \cup L^2 \cup L^3$
and the edge set $K = K^s$ is a 
6-blockade.

 We know that $c(R_1) = 5$ and by the induction hypothesis, $c(R_{s-1}) = 6$ if $s \ge 3$.

Suppose first that $M$ has an edge $ab$ in $R_s - Y$.
Then $ab \in E(R_{s-1})$, and so $M$ contains a subset $M'$ which is a  
blockade
in $Q_{s-1}$ separating $a$ and $b$. 
If $|M'| = 6$, then $6 \le |M'| \le |M|$, and we are done.
Since $ab \in E(R_{s-1})$, the sizes of minimum 
blockades
separating $a$ and $b$ in $Q_{s-1}$ and 
$R_{s-1}$ are the same. So if $s \ge 3$, then by the induction hypothesis, $|M'| = 6$, and we are done.

So let $s = 2$ and $|M'| < 6$. Since $c(R_1) = 5$, $|M'|$ is equal to the size of a minimum  
blockade
 in $R_1$, and so $|M'| = 5$.
Since no minimum 
blockade
 of  $R_1$ is a cut in 
$R_2$, we have $5 = |M'| < |M|$, and so $|M| = 6$.

Now suppose that $M$ has no edge in $R_s - Y$, i.e.,
$M \subseteq E(Y) \cup K$. Then it is easy to check (using the symmetries of $Y$) that every such 
blockade
has six edges, and so $|M| = 6$.
\\[1ex]
${\bf (p3)}$
Finally, we prove $(a3)$.
Let, as in ${\bf (p2)}$, $Y = Y_s$ and $X = X_s$. Obviously, each $Y_i$ is isomorphic to $Y$.
It is easy to prove the following claims.
\\[1ex]
\indent
{\sc Claim} 1. {\em 
$Y - L - e$ has a $\Lambda $-factor
for every edge $e$ in $Y$ and 
for every 3-vertex path $L$ in $Y - e$.}
\\[1ex]
\indent
{\sc Claim} 2.
{\em Let $a$ and $b$ be distinct vertices of degree two in $Y$ that belong to different components of $Y - X$ and let $Y(ab)$ be the graph obtained from $Y$ by adding a new edge  $ab$. Let $S$ be a 3-vertex path in $Y(ab)$ containing edge $ab$. Then $Y(ab) - S - e$ has a 
$\Lambda $-factor
for every edge $e$ in $Y(ab) - E(S)$.}
\\[1ex]
\indent
Let 
$D = E(C_s) - E(\cup \{L_k: k \in \{1, \ldots , 3s\})$.
Obviously, $D$ is a matching and
$R_s = \cup \{Y_i: i \in \{1, \ldots , s\}\} \cup D$.
Now, if $P_i$ is a $\Lambda $-factor of $Y_i$, then clearly
$\cup \{P_i: i \in \{1, \ldots , s\}\}$ is a $\Lambda $-factor of $R_s$.
 Therefore from {\sc Claim} 1 we have:
 \\[1ex]
\indent
{\sc Claim} 3. {\em $R_s - L - e$ has a $\Lambda $-factor for every edge $e$ in $R_s$ and every
3-vertex path $L$ in $R_s - D- e$.}

Thus, it remains to consider a pair $(F, e)$, where $F$
is a 3-vertex path containing an edge in $D$ and 
$e \in E(R_s - F)$.
Since $D$ is a matching, $F$ has exactly one edge
$d$ in $D$. 
By symmetry of $R_s$, we can assume that
$d = z_1^3z_2^1$ and $F = z_1^3z_2^1v$, 
where $z_1^3$ and $z_2^1$ are vertices of degree two in $Y_1$ and $Y_2$, respectively, and
$v \in \{x_2^1, z_2^2\}$.

For $z \in V(C_s)$, let $\alpha (z) = x$ if $xz \in E(R_s)$ and $x \in X$. Obviously, $\alpha $ is a function from 
$V(C_s)$ onto $X$. 
\\[0.7ex]
${\bf (p3.1)}$
Suppose first that $v = x_2^1$.
Let $F_{2i -1} = z_{2i-1}^3z_{2i}^1\alpha (z_{2i}^1)$, and so $F_1 = F$. 
 Let $M = \cup \{F_{2i -1}: \{i \in \{1, \ldots , s\}\}$. 
 By symmetry, we can assume that $e \not \in E(M)$.
 Obviously, $M$ is a $\Lambda $-packing in $R_s$ and
 $Y_i \cap M$ consists of three vertices $a_i$, $b_i$, and $x_i$ and one edge $a_ix_i$, where $a_i$ and $b_i$ have degree two in $Y_i$ and belong to different components of $Y_i - X_i$ and $x_i \in X_i$. Let $S_i = x_ia_ib_i$. Then $S_i$ is a 3-vertex path in $Y_i(a_ib_i)$.
 By {\sc Claim} 2, $Y_i(a_ib_i) - S_i - e_i$ has a $\Lambda $-factor $N_i$ for every edge $e_i$ in 
 $Y_i(a_ib_i) - E(S_i)$.
Let $N = \cup \{N_i: i \in \{1, \ldots , s\}\}$.
Since $V(Y_i \cap M) = V(S_i)$, we have:
$M \cup N$ is a $\Lambda $-factor of $R_s$ containing $F$ and avoiding $e$. 
\\[0.7ex]
${\bf (p3.2)}$
Finally, suppose that $v = z_2^2$.
Let $F_i  = z_{2i-1}^3z_{2i}^1z_{2i}^2$, and so $F = F_1$.
Let $M = \cup \{F_i: i \in \{1, \ldots , s\}\}$.
 By symmetry, we can assume that $e \not \in E(M)$.
As in ${\bf (p3.1)}$, $M$ is a $\Lambda $-packing in $R_s$ and
 $Y_i \cap M$ consists of three vertices $a_i$, $b_i$, and $x_i$ and one edge $a_ix_i$, where $a_i$ and $b_i$ have degree two in $Y_i$ and belong to different components of $Y_i - X_i$ and $x_i \in X_i$. Therefore we are done by the same arguments as in ${\bf (p3.1)}$.
 \ep
\\[3ex]
{\bf Acknowledgment.} 
I am grateful to  the referees for a very careful reading of the paper and  very useful remarks.
Also, many thanks to one of the referees for drawing my attention to recent papers \cite{Kos1,Kos2,Kos3}.


\end{document}